\documentclass[red,11pt,a4paper]{article}

\usepackage{cite}

\usepackage{amsmath}
\usepackage{amscd}
\usepackage{amssymb}
\usepackage[pdftex]{color,graphicx,hyperref}
\usepackage{latexsym}
\usepackage{color}
\usepackage{cases}
\usepackage[normalem]{ulem}

\usepackage[latin1]{inputenc}
\usepackage[T1]{fontenc} %only using latex
\usepackage[english]{babel}
\usepackage{verbatim}
\usepackage{calligra}
\usepackage{enumerate}
\usepackage{dsfont}
\usepackage{amsfonts}
\usepackage{mathtools}
\usepackage{mathrsfs}
\usepackage{bm}

\bf
\usepackage{ika}

\newcommand{\T}{\mathbb{T}}

\def\tD{{\mathbb D}}

\newcommand{\tsl}{\textsl}

\newcommand{\mbb}{\mathbb}

\newcommand{\mc}{\mathcal}

\newcommand{\veps}{\varepsilon}
\newcommand{\eps}{\veps}
\newcommand{\what}{\widehat}
\newcommand{\wtilde}{\widetilde}
\newcommand{\vphi}{\varphi}
\newcommand{\oline}{\overline}
\newcommand{\ra}{\rightarrow}

\newcommand{\g}{\gamma}
\renewcommand{\k}{\kappa}

\newcommand{\de}{\delta}

\newcommand{\lan}{\langle}
\newcommand{\ran}{\rangle}
\newcommand{\ess}{{\rm ess}}
\newcommand{\res}{{\rm res}}

\newcommand{\Q}{\mathbb{Q}}

\newcommand{\N}{\mathbb{N}}
\newcommand{\Z}{\mathbb{Z}}
\renewcommand{\P}{\mathbb{P}}

\renewcommand{\div}{{\rm div}\,}

\newcommand{\curl}{{\rm curl}\,}

\newcommand{\Id}{{\rm Id}\,}
\newcommand{\Supp}{{\rm Supp}\,}

\allowdisplaybreaks

\def\d{\partial}
\def\div{{\rm div}\,}

\begin{document}

\title{Low Mach number limit for degenerate Navier-Stokes equations \\
in presence of strong stratification}

\author{Francesco Fanelli
\thanks{The work of F.F. has been partially supported by the LABEX MILYON (ANR-10-LABX-0070) of Universit\'e de Lyon, within the program ``Investissement d'Avenir''
(ANR-11-IDEX-0007),  and by the projects BORDS (ANR-16-CE40-0027-01), SingFlows (ANR-18-CE40-0027) and CRISIS (ANR-20-CE40-0020-01), all operated by the French National Research Agency (ANR).}
\and Ewelina Zatorska
\thanks{The work of E.Z. was partially supported by the EPSRC Early Career Fellowship no. EP/V000586/1.}
}

\date{\today}

\maketitle

{
\footnotesize
\centerline{$^*\;$Univ. Lyon, Universit\'e Claude Bernard Lyon 1, CNRS UMR 5208, Institut Camille Jordan,}
\centerline{43, Boulevard du 11 novembre 1918 -- F-69622 Villeurbanne, France}
\centerline{\small \texttt{fanelli@math.univ-lyon1.fr}}

\bigbreak
\centerline{$^\dagger\;$Department of Mathematics, Imperial College London}
\centerline{6M14 Huxley Building, South Kensington Campus -- SW7 2AZ, London, UK}
\centerline{\small \texttt{e.zatorska@imperial.ac.uk}}

}

\bigbreak
\bigbreak

\begin{abstract}

In this paper, we investigate the low Mach and low Froude numbers limit for the compressible Navier-Stokes equations with degenerate, density-dependent, viscosity coefficient, in the strong stratification regime. We consider the case of a general pressure law with singular component close to vacuum, and general ill-prepared initial data. We perform our study in the three-dimensional periodic domain. We rigorously justify the convergence to the generalised anelastic approximation, which is used extensively to model atmospheric flows.

\end{abstract}


\medbreak
(secondary).
\bigbreak

\section{Introduction} \label{s:intro}

Flows in the atmosphere are typically characterised by two main features (see \cite{Ped}): first of all, they are weakly compressible, moreover
they undergo the combined effect of a strong stratification (due to the action of gravity) and of a strong Coriolis force (due to the rotation of the Earth, which
is very fast if compared to the space-time scales of the flows).

Neglecting the effects of the Earth rotation, the importance of the other two factors, \tsl{i.e.} weak compressibility and strong stratification, may be assessed
by introducing two physical a-dimensional parameters, the Mach number and the Froude number, respectively. The smaller these parameters are, the more predominant
 weak compressibility and strong stratification become. Thus, as usual in Physics, it is natural to look at the regime where both parameters vanish, to find reduced models for atmospheric flows. They are simpler to deal with than the corresponding primitive system, both from the analytical and numerical point of view. 

When the Mach number and the Froude number go to zero with the same speed, the flow becomes incompressible and stratified at the same rate. Formally, this asymptotic regime was considered already by Ogura and Phillips in \cite{O-P}. The limiting system takes the name of anelastic approximation. The
physical importance of the anelastic approximation  is discussed, for example in 
\cite{Klein-10} in the context of various atmospheric flows, and in \cite{ABRZ} in the context of astrophysics models.

\subsection{The primitive system and the limit dynamics}
In this paper we will give a rigorous derivation of what we call
the generalised anelastic approximation, namely an anelastic approximation with variable viscosity. The
starting system (referred to as the primitive system) is the barotropic Navier-Stokes equations, with bulk viscosity coefficient equal to $0$ and the shear viscosity coefficient proportional to the density of the fluid. In particular, the system strongly degenerates close to vacuum.
This choice of the viscosity coefficients is physically relevant, as viscosity is, in general, hardly expected to be uniform for flows on large scales. Their specific form allows one to exploit a certain mathematical structure of the system, called the BD-entropy (see more details in the discussion below).

Assuming that both the Mach and Froude numbers are equal to a small parameter $\veps>0$, the system of equations reads as follows:
\eq{\label{main}
&\partial_t \vr + \Div (\vr \vu) = 0,\\
&\partial_t (\vr \vu) + \Div (\vr \vu \otimes \vu) + \frac{1}{\ep^2}\Grad p(\vr) -  \nu\Div (\vr\,\tD\vu) = \frac{1}{\ep^2} \vr \Grad G.
}
The unknowns are the mass density $\vr\,=\,\vr(t,x)\,\geq\,0$, and the velocity vector  field $\vu\,=\,\vu(t,x)\in\R^3$. The function $p\,=\,p(\vr)$ denotes
the internal pressure, the constant $\nu>0$ is the viscosity coefficient, and $\tD=\frac12\lr{\Grad+\Grad^t}$ is the symmetric part of the gradient. Finally, $G=G(x)$ is a smooth function
(say $G\in C^3(\Omega)$) describing a scalar external force acting on the flow. $G$ typically encodes the action of gravity, in which case
$G=-gx_3 $, where $g$ is the gravitational acceleration constant.

Due to the present state-of-the-art of the mathematical theory for system \eqref{main}, 
we assume that
the fluid occupies the periodic box in $\R^3$, \tsl{i.e.} we consider the equations  \eqref{main} on the space domain
\begin{equation} \label{eq:domain}
 \Omega\,:=\,\T^3\,.
\end{equation}

The pressure function $p$
is assumed of the following form:
\eq{ \label{eq:pressure}
p(\vr)\,=\,p_e(\vr)+p_c(\vr)\,=\,\dfrac{1}{\g}\,\vr^\gamma-\frac{1}{\kappa}\vr^{-\kappa}\,, \qquad\qquad \gamma>1,\quad \kappa\geq \gamma-2, \quad \k>3.}
The first part is the standard barotropic pressure, while the second part is the so-called ``cold pressure'', because it is most significant in the region of temperatures close to zero.
The constants $1/\g$ and $1/\k$ are just normalisation factors; their presence is useful in some computations.
The restriction on the adiabatic exponent $\g>1$  is somehow classical in the theory of compressible Navier-Stokes equations. The conditions on the exponent $\k$, instead, are of technical nature.

\medbreak

When $\veps\to0^+$ in equations \eqref{main}, we observe a competition between the large size of the pressure term (low Mach number effect), which
tends to drive the flow to incompressibility, and the large size of the forcing term (low Froude number effect), which increases the stratification of the flow.
Due to the choice of scaling, those two terms are in balance in the limit process. Therefore, it is easy to see that,
when $\ep\to 0^+$, $\vr$ will tend, at least formally, to the profile $b=b(x)$ satisfying 
\begin{equation} \label{eq:b}
\nabla p(b)\,=\,b\,\nabla G\,.
\end{equation}
Smoothness of $G(x)$ and strict monotonicity of $p$ imply that there exists a smooth function $b\in C^3(\Omega)$ satisfying \eqref{eq:b}. Monotonicity of $p$ implies convexity of the pressure potential $H$ (defined in \eqref{defH} below) which provides existence of constants $b_*, b^*\,\in\R$ such that
\begin{equation} \label{est:b}
\forall\,x\in\Omega\,,\qquad\qquad 0\,<\,b_*\,\leq b(x)\,\leq\,b^*\,.
\end{equation}
Note that, for $p(\vr)=\frac{\vr^2}{2}$, one gets $G=b$ up to the choice of an irrelevant additive constant, which is the case considered in \cite{BGL}.
On the other hand, if $G(x)\,=\,-gx_3$ is the gravitational potential, it is easy to see that $b=b(x_3)$ verifies \eqref{est:b}.

Since $\vr\approx b$ for $\ep\to0^+$, assuming that we can identify the limits of the non-linear terms appearing in \eqref{main}, 
it is not difficult to check that (formally) the limiting system is an anelastic approximation with variable viscosity coefficient, namely
\eq{ \label{eq:limit}
&\Div(b\,\vc U)\,=\,0 \\[1ex]
&\pt\vc U\,+\,(\vc U\cdot\Grad)\vc U\,+\,\Grad\pi\,-\,b^{-1}\,\nu\,\Div\big(b\,\tD\vc U\big)\,=\,0\,.
}
We refer to system \eqref{eq:limit} as the generalised anelastic approximation.
In the above system, $\pi\,=\,\pi(t,x)$ denotes an unknown scalar function, and the term $\nabla\pi$ plays the role of a Lagrangian multiplier associated with the anelastic constraint
$\div(b\,\vc U)=0$.  The limiting system can be also regarded as the viscous counterpart of the so-called lake equation, whose study was initiated in \cite{L-O-T}.

The goal of this paper is to rigorously justify the above formal derivation in the framework of \emph{global in time finite energy weak solutions} to the primitive system \eqref{main}-\eqref{eq:pressure} for general \emph{ill-prepared initial data}. 

\subsection{An overview of related results}

Due to the physical importance of the anelastic approximation,  its rigorous derivation has been the object of intensive studies in the past years.

In \cite{Masm}, Masmoudi proved the rigorous derivation of the anelastic approximation, starting from the classical barotropic Navier-Stokes system. He considered a bounded domain in $\R^3$, supplemented
with Dirichlet boundary conditions, and the limit was performed for ill-prepared initial data \tsl{via} a compensated compactness argument.
Soon after that, Feireisl, M\'alek, Novotn\'y and Stra\v{s}kraba proved an analogous result on a periodic box, and for pressure laws which are
small variations of the ideal gas law, see \cite{F-M-N-S}. Finally, we refer to the book by Feireisl and Novotn\'y \cite{F-N} for a complete account of the mathematical literature on the low Mach number limit, in the presence of both low and high stratification effects. 
They presented the theory for the full Navier-Stokes-Fourier system and in the framework of global in time finite energy weak solutions.  The literature related to the incompressible limit of compressible fluid equations is of course much more extensive
(see \tsl{e.g.} the pioneering works \cite{Ebin} by Ebin and \cite{K-M_1981}-\cite{K-M_1982} by Klainerman and Majda) and recalling all relevant results  goes far beyond the scope of this introduction.
We thus limit ourselves to quote a couple of recent works.

In \cite{F-K-N-Z}, a variant of the anelastic approximation was derived, starting from a version of the Navier-Stokes-Fourier system with neglected thermal diffusion: the potential temperature is assumed to be just transported by the velocity field. The limit system that is identified in \cite{F-K-N-Z} reads as a coupling of the anelastic approximation system \eqref{eq:limit} with a transport equation for the limiting temperature. The convergence is proven in the infinite slab $\R^2\times\,]0,1[\,$ through a spectral analysis of the singular perturbation operator and an application of the celebrated RAGE theorem from scattering theory. The advantage of that technique, in comparison to the one used in \cite{Masm}, is that it allows to get the compactness of the sequence of velocity fields.

On the side of the incompressible limit (with no stratification effects, though),
another interesting result is  \cite{D-M}, where the the authors deal with weakly compressible viscous fluids
in a critical regularity functional framework. In that paper, weak compressibility is obtained by taking a large bulk viscosity coefficient limit,
instead of the classical low Mach number limit. More recently, in \cite{F_2019}, a similar idea was implemented for fast rotating fluids.

For the degenerate Navier-Stokes system \eqref{main}, as considered in our paper, the relevant results are much more sparse. The first one to mention is
the incompressible limit in a strong stratification regime considered in \cite{BGL} by Bresch, Gisclon and Lin. Their system includes two artificial drag terms
in the momentum equation, in order to improve the available information for the velocity field close to vacuum.
The convergence to the anelastic approximation is proven using the relative energy method, for a special choice of pressure law $p(\vr)=\vr^2/2$  and for well-prepared initial data.

A similar method was used in \cite{B-D_2003} by Bresch and Desjardins and in \cite{J-L-W} by J\"ungel, Lin and Wu, for the low Mach number and low Rossby number limit in a two-dimensional geometry. The external force in these papers is replaced by the Coriolis force (whence the low Rossby number regime) and a capillarity term.
The resulting limiting system is a quasi-geostrophic type equation for the stream function
of the limit velocity field. We also refer to \cite{F_MA}-\cite{F_JMFM} for a generalisation of these results to the $3$-D setting and to the case of general ill-prepared initial data.

\subsection{The content of the paper}

In the context depicted above, our work can be seen as a generalisation of the result from \cite{BGL}, to the case of ill-prepared initial data and of more general pressure laws (and hence, more general external forces $G$). 
We work in the framework of global in time finite energy weak solutions to the primitive system. Their existence, in presence of a cold part of the pressure \eqref{eq:pressure}, has been established in \cite{Zatorska_JDE}, \cite{MuPoZa_DCDS}. The case without this assumption has been completed much more recently in \cite{V-Yu} by Vasseur and Yu.
In all these results, the finite energy condition plays, of course, a major role.
However, the degenerate Navier-Stokes system \eqref{main} possesses also a second energy inequality, usually named \emph{BD entropy inequality}
after Bresch and Desjardins, who investigated this second energy conservation law in \cite{B-D_2003}.

The BD entropy estimate provides a control on the gradient of a certain function of the density, whose exact form depends on the form of the viscosity coefficient. For system \eqref{main}, this function is $\Grad \sqrt{\vr}$. The BD entropy also allows to control the skew-symmetric part $\mbb A=\frac12\lr{\Grad-\Grad^t}$ of the gradient of the velocity. This, when combined with the classical energy, provides the corresponding bound for the full gradient of the velocity.

The classical energy inequality, the BD entropy inequality, and all the bounds that follow, are essential also in the present paper.
As a matter of fact, for any value of parameter $\veps\in\,]0,1]$, we consider a finite energy weak solution $\big(\vr_\ep,\vu_\ep\big)$
to system \eqref{main}, which satisfy both those energy inequalities. However, proving that the BD entropy estimate is satisfied uniformly with respect to $\ep$ requires some effort, especially when a general pressure law is considered: this is one of the first problems solved in our paper.

Having all these estimates satisfied uniformly for a sequence of
finite energy weak solutions $\big(\vr_\ep,\vu_\ep\big)_\ep$, the rest of the proof of the derivation of the generalised anelastic approximation \eqref{eq:limit} boils down to showing that the weak limit $(b, \vc U)$ is indeed a solution to \eqref{eq:limit}. It is well known that passing to the limit in the weak formulation of equations \eqref{main}, especially in its nonlinear parts, is problematic. This is because
the singular terms, along with the ill-prepared initial data, are responsible for fast time oscillations of the solutions (the so-called \emph{acoustic waves}), which may prevent, in the end, the convergence of the nonlinear terms to the expected limit. Showing that this does not happen is the core of the whole proof.

The main concern is the convergence of the convective term in the momentum equation. To that purpose, we use a different technique than the one from \cite{BGL}. Our approach is inspired by the previous works \cite{L-M} and \cite{Masm} on the incompressible limit for the classical barotropic Navier-Stokes system, and is based on a compensated compactness argument. More precisely, we first regularise the primitive equations, which we recast in the form of a wave system. After that, we exploit two pieces of information coming from the wave system. First of all, we may deduce the compactness of the rotational part of the velocity fields. On the other hand, by direct but elaborated algebraic manipulations, we may infer that the interaction of the potential part of the velocity fields in the convective term gives rise to small quantities, which tend to vanish when $\veps\ra0^+$.
It is worth to point that this argument is robust enough to deal with other variants of the system \eqref{main}. For instance, we could trade
the cold component of the pressure function, which basically provides us with some integrability properties for $\Grad\vu$, for a turbulent drag term $\vr|\vu|\vu$, which would give a better integrability of  the momentum $\vc V\,:=\,\vr\vu$. In that case, most of the steps are the same, 
although the derivation of essential estimates becomes significantly more laborious. The only problem, and the breaking point, arises  when one wants to pass to the limit in this artificial drag term. This term turns out to be even more non-linear than the convective term, because of the presence of the norm $|\vu|$ of the vector $\vu$. It is not clear how to bypass this difficulty in our framework, and so, the problem remains open.

\medbreak
We conclude with a short outline of the paper. In the next section, we collect our main hypotheses on the initial data, we give the definition of finite energy weak solutions, and we state our main result. In Section \ref{s:unif-bounds}, we deduce, from the energy inequality and the BD entropy inequality, a long list of uniform bounds
for the family of solutions $\big(\vr_\ep,\vu_\ep\big)_\ep$ we consider. That part of the study is rather delicate, due to the degeneracy of the system close to vacuum.
In Section \ref{s:singular}, we use the previous uniform bounds to extract a weakly convergent subsequence, and to derive first basic properties on its weak limit point. At this stage we reformulate the primitive equations into the wave system, which describes the propagation of the  acoustic waves. In Section \ref{s:convergence} we rigorously perform the convergence
in the weak formulation of equations \eqref{main}, and conclude the derivation of the anelastic approximation \eqref{eq:limit}. For the convenience of the reader, we collect
some tools from Fourier analysis which we need in our study in the Appendix at the end of the paper.

\medbreak

\section{Statement of the main result} \label{s:system}
In this section, we first introduce our assumptions on the initial data, then we define the notion of finite energy weak solutions to system \eqref{main}-\eqref{eq:pressure}, and finally
we formulate our main theorem.

\paragraph{Initial data.}
Problem \eqref{main}- \eqref{eq:pressure} is supplemented by general \emph{ill-prepared} initial data. Namely, for any small parameter $\veps\in\,]0,1]$ fixed, we pick initial data
\eq{\label{initial_data}
(\vr,\vu)|_{t=0}=
\big(\vr_{0,\ep},\vu_{0,\ep}\big)}
satisfying the following conditions:
\begin{itemize}
\item[(i)] the initial densities $\vr_{0,\ep}\geq0$ are assumed to be small perturbations of the static state $b$, defined by \eqref{eq:b}: more precisely, we assume\footnote{Here and throughout
this paper, we make use of the following notation: given a normed space $X$ and a sequence of functions $\big(f_\veps\big)_\veps$ all belonging to $X$,
we write $\big(f_\veps\big)_\veps\,\subset\,X$ implicitly meaning that the sequence is also \emph{bounded} in $X$.} that
\[
\vr_{0,\ep}\,=\,b\,+\,\ep\,\phi_{0,\ep}\,,\qquad\qquad \mbox{ with }\qquad
\big(\phi_{0,\ep}\big)_\ep\subset L^\infty (\Omega)\quad \mbox{ and } \quad \lr{\Grad\ln\frac{\vr_{0,\ep}}{b}}_{\ep}\subset L^\infty(\Omega)\,;
\]
\item[(ii)] the initial velocity fields $\vu_{0,\ep}$ are such that $\big(\vu_{0,\ep}\big)_\ep\subset L^\infty (\Omega)$.
\end{itemize}

Thus, up to extraction of a subsequence, not relabeled here, we may suppose that
\eq{\label{data_limit}
\phi_{0,\ep}\,\to\,\phi_{0}\qquad \mbox{ and } \qquad \mbox \vu_{0,\ep}\to\vu_{0}\qquad\qquad \mbox{ weakly-$*$ \ in }\ L^\infty(\Omega)\,.
}

\paragraph{Energy functionals.}
Next, we need to introduce various  energy functionals. The internal energy function (sometimes called \emph{pressure potential}) is defined by the ODE
\[
\vr\,H'(\vr)\,-\,H(\vr)\,=\,p(\vr)\,,
\]
which implies in particular that
\[
H''(\vr)\,=\,\frac{p'(\vr)}{\vr}\,.
\]
Notice that $H$ is defined up to the sum of an affine function. Here, we fix the classical choice
\eq{\label{defH}
H(\vr)\,=\,\vr\int_{1}^\vr\frac{p(z)}{z^2}\,{\rm d}z\,=\,\frac{1}{\g\,(\g-1)}\,\vr\,\big(\vr^{\g-1}\,-\,1\big)\,
+\,\frac{1}{\kappa\,(\kappa+1)}\,\vr\,\big(\vr^{-\kappa-1}\,-\,1\big).
}
We now denote
\begin{align}
\mc E\Big(\vr,\vu\,\Big|\,b\Big)\,&:=\,\frac{1}{2}\,\intO{\vr\,|\vu|^2}\,+\,\frac{1}{\ep^2}\,\intO{\Big(H(\vr)\,-\,H(b)\,-\,H'(b)\,(\vr-b)\Big)}
\label{def:E} \\[1ex]
\mc F\Big(\vr,\vu\,\Big|\,b\Big)\,&:=\,
\intO{\vr\,\left|\vu\,+\,\nu\Grad\ln\frac{\vr}{b}\right|^2} \label{def:F}
\end{align}
to be the \emph{classical energy} and the \emph{BD entropy} functions.
We also set ${\cal E}\Big(\vr,\vu\,\Big|\,b\Big)(T)\,:=\,{\cal E}\Big(\vr(T),\vu(T)\,\Big|\,b\Big)$ and ${\cal E}\Big(\vr,\vu\,\Big|\,b\Big)(0)\,:=\,{\cal E}\Big(\vr_{0},\vu_{0}\,\Big|\,b\Big)$, 
and similarly for the function $\mc F$.

\paragraph{Weak solutions to the primitive system.}
After this preparation, we are ready to give the definition of weak solutions to system \eqref{main}-\eqref{eq:pressure} which are relevant for us.

\begin{df}\label{df:main}
Let $\big(\vr_0,\vu_0\big)$ be such that $\,{\cal E}\Big(\vr_{0},\vu_{0}\,\Big|\,b\Big)+{\cal F}\Big(\vr_{0},\vu_{0}\,\Big|\,b\Big)<+\infty$. \\
We say that the couple $(\vr,\vu)$ is a \emph{finite energy weak solution} of \eqref{main}-\eqref{eq:pressure} in $[0,T[\,\times \Omega$, with the initial datum $\big(\vr_0,\vu_0\big)$,
provided the following conditions are satisfied:
\begin{enumerate}[(1)]
\item $\vr\geq0$ almost everywhere, with $\vr\in L^\infty\big([0,T[\,;L^\gamma(\Omega)\big)$ and $\vr^{-1}\in L^\infty\big([0,T[\,;L^\kappa(\Omega)\big)$,
$\Grad\sqrt{\vr}\in L^{\infty}\big([0,T[\,;L^2(\Omega)\big)$ and $\sqrt{\frac{p'(\vr)}{\vr}}\Grad\vr\in L^2\big([0,T[\,; L^2(\Omega)\big)$;
\item $\sqrt{\vr}\vu\in L^\infty\big([0,T[\,; L^2(\Omega)\big)$ and $\sqrt{\vr}\Grad\vu\in L^2\big([0,T[\,; L^2(\Omega)\big)$;
\item the equations of system \eqref{main} are satisfied in the sense of distributions: more precisely, we have
\eq{\label{weak:cont}
\intO{\vr_0\xi(0)}+\intTO{\Big(\vr\pt\xi\,+\,\vr\vu\cdot\Grad\xi\Big)}\,=\,0}
for any test function $\xi\in \mc D\big([0,T[\,\times\Omega\big)$, and
\eq{\label{weak:mom}
&\intO{\vr_0\vu_0\cdot\psi(0)}+\intTO{\Big(\vr\vu\cdot\pt\psi+\vr\vu\otimes\vu:\Grad\psi\Big)}\\
&\quad+\frac{1}{\ep^2}\intTO{p(\vr)\Div\psi}+\frac{1}{\ep^2}\intTO{\vr\Grad G\cdot\psi}-\nu\intTO{\vr\,\tD\vu:\Grad\psi}\,=\,0
}
for any test function $\psi\in \mc D\big([0,T[\,\times\Omega;\R^3\big)$;
\item for almost every $t\in[0,T[\,$, the following energy inequalities hold true:
\eq{ \label{est:en-BD-weak}
&\mc E \Big(\vr,\vu\,\Big|\,b\Big)(t)+\nu\int_0^t\!\!\intO{\vr|\tD\vu|^2}\,{\rm d}s\leq \mc E \Big(\vr_0,\vu_0\,\Big|\,b\Big)\,,\\
&\mc F \Big(\vr,\vu\,\Big|\,b\Big)(t)+\,\frac{\nu}{\veps^2}\int_0^t\!\!\intO{b^2\,\frac{p'(\vr)}{\vr}\,\left|\Grad\lr{\frac{\vr}{b}}\right|^2}\,{\rm d}s+\nu\int_0^t\!\!\intO{\vr|\mbb A\vu|^2}\,{\rm d}s\leq C_0\,e^{C_0(1+T)}\,,
}
where the constant $C_0>0$ may depend on $\mc E \Big(\vr_0,\vu_0\,\Big|\,b\Big)$ and $\mc F \Big(\vr_0,\vu_0\,\Big|\,b\Big)$ but is independent of $\ep$.
\end{enumerate}

The solution $(\vr,\vu)$ is said \emph{global in time} if the previous properties hold true for any $T>0$.
\end{df}

For any $\veps\in\,]0,1]$ fixed, the existence of global in time finite energy weak solutions to system \eqref{main} in the sense of previous definition
was proven in \cite{Zatorska_JDE} and \cite{MuPoZa_DCDS}, in the case $G=0$ (corresponding to $b=const.$). The argument of those papers apply in a fairly direct way also to the case 
considered in this paper, where $G\neq 0$ and $b$ is non-constant: we explain in the next section how to modify the estimates of \cite{Zatorska_JDE}-\cite{MuPoZa_DCDS} in order to include the force.

\paragraph{Main result.}
Before stating the main result of this paper, we need some additional tools and notation. Following \cite{Masm}-\cite{F-M-N-S} (see also \cite{Lions_1}), we introduce the
twisted Leray-Helmholtz projector $\P_b$, related to the smooth function $b$ satisfying \eqref{est:b}, as follows: for any smooth vector field $\vc v$ on $\Omega$, we write
\[
\vc v\,=\,\P_b[\vc v]\,+\,b\,\Grad\Psi\,,
\]
where $\Psi$ is the unique solution to the Neumann problem
\[
\Div\big(b\,\Grad\Psi\big)\,=\,\Div\vc v\quad\mbox{ in }\ \Omega\,,\qquad\qquad \intO{\Psi}\,=\,0\,.
\]
Remark that $\P_b[\vc v]$ and $\Q_b[\vc v]\,:=\,b\,\Grad\Psi$ are orthogonal in the weighted Hilbert space $L^2_b(\Omega;\R^3)$, which is defined as the space of functions
$f:\Omega\longrightarrow \R^3$ which are $L^2$-summable with respect to the measure $\frac1b\,\dx$.

Similarly to the case of the classical Leray-Helmholtz projector $\P\,=\,\P_1$ and its $L^2$-orthognal projector $\Q\,=\,\Q_1$, it is possible to prove that both
$\P_b$ and $\Q_b$ are bounded continuous functionals on $L^p(\Omega;\R^3)$, for any $1<p<+\infty$.

\medbreak
We can now state the main result of this paper, which is contained in the following theorem.
\begin{thm}\label{thm:main}
Let $ \gamma>1$ and $\kappa\geq \gamma-2,\; \k>3$ in \eqref{eq:pressure}.\\
Let $\big(\vr_{0,\veps},\vu_{0,\veps}\big)_\veps$ be a family of initial data satisfying hypotheses (i)-(ii) fixed above, so in particular the condition
\begin{equation} \label{ub:data}
\sup_{\veps\in\,]0,1]}\mc E\Big(\vr_{0,\veps},\vc \vu_{0,\veps}\,\Big|\,b\Big)\,+\,
\sup_{\veps\in\,]0,1]}\mc F\Big(\vr_{0,\veps},\vc \vu_{0,\veps}\,\Big|\,b\Big)\,<\,+\infty
\end{equation}
holds. Define the couple $(\phi_0,\vu_0)$ as in \eqref{data_limit}. \\
Let $\big(\vre,\vue\big)_\ep$ be a family of global in time weak solutions to system \eqref{main}-\eqref{eq:pressure}, in the sense of Definition \ref{df:main},
corresponding to the previous initial data.
Define the scalar quantity $\phi_\ep:=\frac{\vre-b}{\ep}$. 

Then, there exists a couple of functions $\big(\phi,\vc U\big)$ such that, passing to a suitable subsequence as the case may be, in the limit $\ep\to 0$ one has
\begin{align*}
&\vr_\ep\,\to\, b\qquad \mbox{ strongly in }\ L^\infty_{\rm loc}\big(\R_+;L^p(\Omega)\big)\,, \quad \mbox{ for any }\ p\in[1,3[\,,\\[1ex]
&\phi_\ep\,\to\,\phi\qquad \mbox{ weakly in}\ L^2_{\rm loc}\big(\R_+; W^{1,\min\{\gamma,2\}}(\Omega)\big)\,,\\[1ex]
&\vue\,\to\,\vc U\qquad \mbox{ weakly-$*$ in }\ L^\infty\big(\R_+;L^{p_1}(\Omega)\big)\,\cap\,L^2_{\rm loc}\big(\R_+;W^{1,p_1}(\Omega)\big),\quad \mbox{ where }\
p_1\,:=\,\frac{2\k}{\k+1}\,.
\end{align*}
In addition, $\phi\,=\,\phi(b)$ is a function of the static profile $b$, while $\vc U$ is a solution of the target system \eqref{eq:limit} in the weak sense,
related to the initial datum $\vc U|_{t=0}\,=\,\vc U_0\,:=\,\frac{1}{b}\,\P_b[b\,\vu_0]$, \tsl{i.e.} one has
$\Div(b\vc U)=0$ almost everywhere in $\R_+\times\Omega$ and
\eq{\label{weak:mom_lim}
\intO{b\,\vu_0\cdot\psi(0)}+\intTO{\Big(b\, \vc U\cdot\pt\psi\,+\,b\,\vc U\otimes\vc U:\Grad\psi\Big)}-\nu\intTO{b\,\Grad\vc U:\Grad\psi}\,=\,0
}
for any $T>0$ and any test function $\psi\in \mc D\big([0,T[\,\times\Omega;\R^3\big)$ such that $\Div(b\psi)=0$.
\end{thm}

\begin{rmk}
Note that the initial condition equals
\eqh{\intO{b\,\vu_0\cdot\zeta}=\intO{\P_b[b\,\vu_0]\cdot\zeta}=\intO{b\,\vc U_0\cdot\zeta}
}
for any test function $\zeta\in \mc D\big(\Omega;\R^3\big)$ such that $\Div(b\zeta)=0$.
\end{rmk}

\section{\textsl{A priori} estimates} \label{s:unif-bounds}

In this section, we derive \emph{uniform bounds} for the family of weak solutions $\big(\vr_\veps,\vu_\veps\big)_\veps$
to the original Navier-Stokes system \eqref{main}. The main tools for this are the classical \emph{energy inequality} and the so-called \emph{BD entropy estimate}.

Here and everywhere in the text, we adopt the following notation: given a Banach space $X$ and any $p\in[1,+\infty]$, we set $L^p_T(X)\,:=\,L^p([0,T];X)$;
in the case $T=+\infty$, instead, we explicitly write $L^p(\R_+;X)$. When convenient, we will use also the notation $L^p_{\rm loc}(\R_+;X)\,:=\,\bigcap_{T>0}L^p_T(X)$.

\subsection{Bounds coming from the energy inequality} \label{ss:energy}

The energy inequality for $\big(\vr_\veps,\vu_\veps\big)$, which is satisfied by assumption, reads as follows: for almost any time $T>0$, we have
\begin{align} \label{est:energy}
{\cal E}\Big(\vr_\veps,\vu_\veps\,\Big|\,b\Big)(T)
+\nu\intTO{\vr_\veps\,|\tD\vu_\veps|^2}\,\leq\,{\cal E}\Big(\vr_{0,\veps},\vc m_{0,\veps}\,\Big|\,b\Big)\,,
\end{align}
where the function $\mc E\Big(\vr,\vu\,\Big|\,b\Big)$ has been defined in \eqref{def:E}
and we recall that we have set ${\cal E}\Big(\vr_\veps,\vu_\veps\,\Big|\,b\Big)(T)\,:=\,{\cal E}\Big(\vr_\veps(T),\vu_\veps(T)\,\Big|\,b\Big)$.
Notice that, due to the cold pressure, at any value of $\veps\in\,]0,1]$ fixed, the velocity field $\vu_\veps$ is well-defined , thus the previous notation makes sense.

\medbreak
From the energy inequality \eqref{est:energy}, we now derive first uniform bounds for the family $\big(\vr_\veps,\vu_\veps\big)_\veps$.
In fact, owing to our assumptions on the initial data, and in particular to \eqref{ub:data}, the right-hand side of \eqref{est:energy} is \emph{uniformly bounded}: specifically, we have
\[%\begin{equation} \label{ub:data}
\sup_{\veps\in\,]0,1]}\mc E\Big(\vr_{0,\veps},\vc \vu_{0,\veps}\,\Big|\,b\Big)\,<\,+\infty\,.
\]%\end{equation}
Then, it is easy to deduce the following uniform bounds:
\begin{align}
\big(\sqrt{\vr_\veps}\,\vu_\veps\big)_\veps\,&\subset\,L^\infty\big(\R_+;L^2(\Omega)\big)\,, \label{ub:kinetic} \\
\big(\sqrt{\vr_\veps}\,\tD \vu_\veps\big)_\veps\,&\subset\,L^2\big(\R_+;L^2(\Omega)\big)\,. \label{ub:Du_zero} 
\end{align}

Let us now focus on the density functions. To begin with, following the approach of \cite{F-N}, it is convenient to decompose any function $h$ into its \emph{essential} and \emph{residual} parts.
Thus, for almost every time $t>0$ and all $\veps\in\,]0,1]$, we introduce the sets
$$
\Omega_\ess^\veps(t)\,:=\,\left\{x\in\Omega\;\Big|\quad \frac{b_*}{2}\,\leq\,\vr_\veps(t,x)\,\leq\,2\,b^*\right\}\,,\qquad\Omega^\veps_\res(t)\,:=\,\Omega\setminus\Omega^\veps_\ess(t)\,,
$$
where the constants $b_*$ and $b^*$ have been defined in \eqref{est:b}
Then, given a function $h$, we can write
$$
h\,=\,\left[h\right]_\ess\,+\,\left[h\right]_\res\,,\qquad\qquad\mbox{ where }\qquad \left[h\right]_\ess\,:=\,h\,\mathds{1}_{\Omega_\ess^\veps(t)}\,.
$$
Here above, $\mathds{1}_A$ denotes the characteristic function of a set $A\subset\Omega$.

For later use, it is convenient to divide the residual set $\Omega^\veps_\res(t)$ further: we define
\[
\Omega^\veps_{\res,B}(t)\,:=\,\left\{x\in\Omega^\veps_\res(t)\;\Big|\quad 0\,\leq\,\vr_\veps(t,x)\,<\,\frac{b_*}{2}\right\}\,,\qquad
\Omega^\veps_{\res,UB}(t)\,:=\,\left\{x\in\Omega^\veps_\res(t)\;\Big|\quad \vr_\veps(t,x)\,>\,2\,b^*\right\}
\]
as the regions where $\vr_\veps$, respectively, stays bounded and may become unbounded.

\medbreak
After this preparation, let us come back to \eqref{est:energy} and derive uniform bounds for the family $\big(\vr_\veps\big)_\veps$.
In the essential set, we can perform  Taylor's expansion of the function $H$; we thus get
\[
\Big[H(\vr)\,-\,H(b)\,-\,H'(b)\,(\vr-b)\Big]_\ess\,\geq\,c\,\big|\vr_\veps\,-\,b\big|^2\,\mathds{1}_{\Omega^\veps_\ess(t)}\,,
\]
which implies that
\begin{equation} \label{ub:dens-ess}
 \sup_{t\in\R_+}\left\|\frac{1}{\veps}\;\left[\vr_\veps-b\right]_\ess\right\|_{L^2(\Omega)}\,\leq\,C \,.
\end{equation}
On the other hand, using the convexity of the function $H$, the fact that $\left|\big[\vr_\veps-b\big]_\res\right|\geq b_*/2$ and equation \eqref{defH},
we discover (see \tsl{e.g.} \cite{F-N-S} for details) the following bounds on the residual set:
\begin{align}
 \sup_{t\in\R_+}\left\|\left[\vr_\veps\right]_\res\right\|^\g_{L^\g(\Omega)}\,
 + \sup_{t\in\R_+}\left\|\left[\vr_\veps^{-1}\right]_\res\right\|^\kappa_{L^\kappa(\Omega)}\,
 +\,\sup_{t\in\R_+}\left\|\left[1\right]_\res\right\|_{L^1(\Omega)}\,\leq\,C\,\veps^{2}\,.\label{ub:dens-res}
\end{align}
The previous estimate immediately implies that
\begin{equation} \label{ub:residual-set}
\sup_{t\in\R_+}\mc L\big(\Omega^\veps_\res(t)\big)\,\leq\,C\,\veps^2\,,
\end{equation}
where we have denoted by $\mc L(A)$ the Lebesgue measure of a set $A\subset\Omega$.

At this point, let us define the quantity
\[
 \phi_\veps\,=\,\frac{1}{\veps}\,\big(\vr_\veps\,-\,b\big)\,.
\]
From the uniform bound \eqref{ub:dens-ess}, we may deduce that
\begin{equation} \label{ub:phi_ess}
\left\|\big[\phi_\veps\big]_\ess\right\|_{L^\infty(\R_+;L^2)}\,\leq\,C\,.
\end{equation}
On the other hand, using \eqref{ub:dens-res} we can compute that for any $p\leq \gamma$, we have
\[
\int_\Omega\left|\big[\phi_\veps\big]_\res\right|^p\,dx\,\leq\,\frac{C}{\veps^{p}}\left(\int_\Omega\left|\big[\vr_\veps\big]_\res\right|^p\,dx\,+\,\int_\Omega[1]_\res\,dx\right)\,\leq\,
C\,\veps^{2-p}\,.
\]
Thus, we finally infer that
\begin{equation} \label{ub:phi_res}
\forall\;1\leq p \leq \g\,,\qquad\qquad
\left\|\big[\phi_\veps\big]_\res\right\|_{L^\infty(\R_+;L^p)}\,\leq\,C\,\veps^{(2-p)/p}\,.
\end{equation}
Of course, this estimate will be useful only in the case when $p$ satisfies the additional restriction $p\leq2$.

\subsection{The Bresch-Desjardins estimate} \label{ss:BD}

As it is apparent from the bounds of the previous subsection, the difficulty with system \eqref{main} is that we lose any control on the velocity fields $\vu_\veps$ and their gradient
$\nabla\vu_\veps$ close to vacuum, specifically in the region $\Omega^\veps_{\res,B}$.
The cold pressure term $p_c$ is of great help in order to bypass that difficulty.

However, the cold pressure term alone is not strong enough to give us all the pieces of information we need to pass to the limit.
On the other hand, system \eqref{main} possesses a nice underlying structure, as evidenced for the first time by Bresch and Desjardins (see \tsl{e.g.} \cite{B-D-L}, \cite{B-D_2003}). By taking advantage of
that structure, it is possible to derive, \textsl{via} the so-called \emph{BD entropy estimates}, some uniform estimates on the gradient of the density functions $\vr_\veps$. This is the goal
of the next lemma.

A bound coming from BD entropy estimates has been required in the definition of weak solutions, see point (4) in Definition \ref{df:main}.
Here we show that such a bound holds \emph{uniformly} with respect to the small parameter $\veps\in\,]0,1]$.
Note that the result in Lemma \ref{l:BD} is stated for smooth solutions to the Navier-Stokes system \eqref{main}. This is solely to justify the manipulations required to derive the inequality. Once the inequality is proven for the smooth solutions, it is possible to deduce that
it is inherited also by the finite energy weak solutions considered in this paper (see \cite{Zatorska_JDE} and\cite{MuPoZa_DCDS} for details).

In the next statement, we resort to the notation introduced in \eqref{def:F}, and we recall that we denote by $\mbb A\vu\,=\,\big(\nabla\vu\,-\,\nabla^t\vu\big)/2$
the skew-symmetric part of the Jacobian matrix of the vector field $\vu$.

\begin{lemma} \label{l:BD}
Let $(\vre,\vue)$ be the smooth solution to \eqref{main}-\eqref{eq:pressure}. Then we have 
the inequality
\eq{\label{BD_estimate}
&\hspace{-1cm}\sup_{t\in\,]0,T[}\left(\mc F\Big(\vr_\veps,\vu_\veps\,\Big|\,b\Big)(t)\,+\,
\frac{1}{\ep^2}\intOB{H(\vr_\veps)\,-\,H(b)\,-\,H'(b)(\vr_\veps-b)}\right)\\
&\qquad\qquad\qquad
+\,\frac{\nu}{\veps^2}\intTO{b^2\,\frac{p'(\vr_\veps)}{\vr_\veps}\,\left|\Grad\lr{\frac{\vr_\veps}{b}}\right|^2}\,+\,\nu\intTO{\vr_\veps\,|\mbb A\vu_\veps|^2}\,\leq\,C\,,
}
where the constant $C>0$ depends only on the initial data and on $T$. In particular, the previous bound is uniform with respect to $\ep\in\,]0,1]$.
\end{lemma}

\pf The proof of this estimate follows closely \cite{BGL}, with the only modifications associated with more general forms of the pressure and of the force.
We proceed in several steps, assuming that $\vr_\ep,\vu_\ep$ are smooth enough to justify all the computations. 
For notational simplicity, in what follows we write $(\vr,\vu)$ instead of $\big(\vr_\veps,\vu_\veps\big)$.

{\it{Step 1}.}  From Lemma 5.1 in \cite{BGL} it follows that for sufficiently smooth solutions
of the continuity equation in \eqref{main}, the following equality holds true
\begin{align*}
&\frac12\Dt \intO{\vr\left|\Grad\ln\frac\vr b\right|^2}
+\intO{\vr\Grad\vu\cdot\Grad\ln\frac\vr b\Grad\ln\frac\vr b} \\
&\qquad\qquad+\intO{\vr\Grad\vu\cdot\Grad\ln b\Grad\ln \frac\vr b}+\intO{\vr\vu\cdot\Grad\Grad\ln b\Grad\ln\frac\vr b}
+\intO{\vr\Grad\Div\vu\Grad\ln\frac\vr b}=0.
\end{align*}

{\it{Step 2}.} In this step, one multiplies the momentum equation by $\nu\Grad\ln\frac\vr b$, we get
\[
\nu\intO{\vr\lr{\pt\vu+\vu\cdot\Grad\vu}\cdot\Grad\ln\frac\vr b}
+\nu^2\intO{\vr\tD\vu:\Grad^2\ln\frac\vr b}+\frac{\nu}{\ep^2}\intO{\lr{\Grad p(\vr)-\vr\Grad G}\cdot\Grad\ln\frac\vr b}=0.
\]
Let us now rewrite each term from the above expression. First, note that the transport term gives
\[
\nu\intO{\vr\lr{\pt\vu+\vu\cdot\Grad\vu}\cdot\Grad\ln\frac\vr b}=\nu\Dt\intO{b\vu\cdot\Grad\frac\vr b}-\nu\intO{\vr\Grad\vu:\Grad^t\vu}+\nu\intO{\vr\vu\cdot\Grad(\Grad\ln b)\vu}.
\]
For the diffusion term, after noticing that $\nabla^2$ is always symmetric, we have
\begin{align*}
\nu^2\intO{\vr\tD\vu:\Grad^2\ln\frac\vr b}&=\nu^2\intO{\vr\Grad\vu:\Grad^2\ln\frac\vr b}\\
&=\nu^2\intO{b\Grad\vu:\Grad^2\frac\vr b}
-\nu^2\intO{\vr\Grad\vu\cdot\Grad\ln\frac\vr b \Grad\ln \frac\vr b}\\
&=-\nu^2\intO{\vr\Grad\vu\cdot\Grad\ln b\Grad\ln \frac\vr b}-\nu^2\intO{\vr\Grad\Div\vu\Grad\ln\frac\vr b}\\
&\qquad\qquad\qquad-\nu^2\intO{\vr\Grad\vu\cdot\Grad\ln\frac\vr b \Grad\ln \frac\vr b}.
\end{align*}
Finally, for the pressure and force term, using \eqref{eq:b}, we obtain
\[
\nu\intO{\frac{\Grad p(\vr)-\vr\Grad G}{\ep^2}\cdot\Grad\ln\frac\vr b}=\frac\nu{\ep^2}\intO{b^2\frac{p'(\vr)}{\vr}\left|\Grad\frac\vr b\right|^2}+
\frac\nu{\ep^2}\intO{\big(p'(\vr)-p'(b)\big)\Grad b\cdot\Grad\frac\vr b}.
\]

{\it{Step 3}.} Now, we sum up the equalities from the previous steps along with the energy estimate.
In our case, after setting $H(\vr;b)\,:=\,H(\vr)-H'(b)(\vr-b)-H(b)$, the statement of Lemma 5.2 from \cite{BGL} gives
\eq{\label{BD_1}
&\Dt\intOB{%\vr\left|\vu+\nu\Grad\ln\frac{\vr}{b}\right|^2
\mc F\Big(\vr,\vu\,\Big|\,b\Big)+\frac{1}{\ep^2}\,H(\vr;b)}
+\frac{\nu}{\veps^2}\intO{b^2\frac{p'(\vr)}{\vr}\left|\Grad\lr{\frac{\vr}{b}}\right|^2}+\nu\intO{\vr|\tD\vu|^2}\\
&=\nu\intO{\vr\Grad\vu:\Grad^{t}\vu}-\nu\intO{\vr\vu\cdot(\Grad\Grad\ln b)\vu}
-\nu^2\intO{\vr\vu\Grad\Grad\ln b\Grad \ln \frac{\vr}{b}}\\
&\quad-\frac\nu{\ep^2}\intO{\big(p'(\vr)-p'(b)\big)\Grad b\cdot\Grad\frac\vr b}\,.
%-\alpha\intO{\vr|\vu|\vu\cdot\Grad\ln\frac{\vr}{b}}
}
Now notice that
\[
\nu\intO{\vr|\tD\vu|^2}-\nu\intO{\vr\Grad\vu:\Grad^{t}\vu}=\nu\intO{\vr|\mathbb{A}\vu|^2},
\]
and so, we finally get
\eq{\label{BD}
&\Dt\intOB{%\vr\left|\vu+\nu\Grad\ln\frac{\vr}{b}\right|^2
\mc F\Big(\vr,\vu\,\Big|\,b\Big)+\frac{1}{\ep^2}\,H(\vr;b)}
+\frac{\nu}{\veps^2}\intO{b^2\frac{p'(\vr)}{\vr}\left|\Grad\lr{\frac{\vr}{b}}\right|^2}+\nu\intO{\vr|\mathbb{A}\vu|^2}\\
&=-\nu\intO{\vr\vu\cdot(\Grad\Grad\ln b)\vu}
-\nu^2\intO{\vr\vu\Grad\Grad\ln b\Grad \ln \frac{\vr}{b}}-\frac\nu{\ep^2}\intO{\big(p'(\vr)-p'(b)\big)\Grad b\cdot\Grad\frac\vr b}\\
%-\alpha\intO{\vr|\vu|\vu\cdot\Grad\ln\frac{\vr}{b}}
&=\sum_{i=1}^3J_i.
}

{\it{Step 4}.} We have to control the terms $J_1,J_2, J_3$ appearing in the right-hand side of \eqref{BD}. To begin with,
we can easily estimate
\begin{align*}
\intT{J_1}&=-\nu\intTO{\vr\vu\cdot(\Grad\Grad\ln b)\vu}\leq\|\Grad^2\ln b\|_{L^\infty_T(L^\infty)}\intTO{\vr|\vu|^2}\leq C\,.
\end{align*}
For $J_2$, instead, we get
\begin{align*}
\intT{J_2}&=-\nu^2\intTO{\vr\vu\Grad\Grad\ln b\Grad \ln \frac{h}{b}}\,\leq\,  \nu^2\|\Grad^2\ln b\|_{L^\infty_T(L^\infty)}
%\intTO{\vr\left|\vu+\nu\Grad\ln\frac{\vr}{b}\right|^2}\,,
\int^T_0\mc F\Big(\vr,\vu\,\Big|\,b\Big)\,\dt\,,
\end{align*}
hence $J_2$ can be controlled by means of a Gr\"onwall argument.\\
Finally, we have to deal with $J_3$. This estimate is a bit more involved than the previous ones, as this term depends on the pressure. We have to distinguish some different cases.

\underline{Case 1: integral over $\Omega^\veps_\ess$}. We start by bounding the part of $J_3$ which is restricted to the essential set.
For this, we use the Taylor expansion $p'(\vr)=p'(b)+p''(z)(\vr-b)$, for some $z$ between $\vr$ and $b$, and the fact that $\vr$ is bounded in $\Omega_\ess$ to write
\begin{align*}
\left|\frac\nu{\ep^2}\int_{\Omega_\ess^\veps}{\big(p'(\vr)-p'(b)\big)\Grad b\cdot\Grad\frac\vr b}\,\dx\right|\,&\lesssim\,
\frac{\nu}{\veps^2}\int_{\Omega_\ess^\veps}|\vr-b|\,|\nabla b| \left|\sqrt{\frac{p'(\vr)}{\vr}}\,b\,\nabla\left(\frac{\vr}{b}\right)\right|\,\dx \\
&\lesssim\,C_\de\,\left\|\big[\phi\big]_\ess\right\|^2_{L^2}\,+\,\delta\,\frac{\nu}{\veps^2}\intO{b^2\,\frac{p'(\vr)}{\vr}\,\left|\nabla\left(\frac{\vr}{b}\right)\right|^2}\,,
\end{align*}
where we have also applied the Young inequality in the last step. Here, $\delta>0$ can be taken arbitrarily small, to the price of increasing the value of $C_\delta$. Thus,
taking $\delta$ small enough, after integration in time, we can absorbe the last term of the previous estimate in the left-hand side of \eqref{BD},
while the other term is bounded by a uniform constant $C$ times $T$, in view of \eqref{ub:phi_ess}.

\underline{Case 2: integral over $\Omega_\res^\veps$}.
In the residual set, it is convenient to split $J_3$ into two pieces, as follows:
\begin{align}
\frac\nu{\ep^2}\int_{\Omega_\res^\veps}{\big(p'(\vr)-p'(b)\big)\Grad b\cdot\Grad\frac\vr b}\,\dx\,&=\,
\frac\nu{\ep^2}\int_{\Omega_\res^\veps}{p'(\vr)\Grad b\cdot\Grad\frac\vr b}\,\dx\,-\,\frac\nu{\ep^2}\int_{\Omega_\res^\veps}{p'(b)\Grad b\cdot\Grad\frac\vr b}\,\dx \label{est:J4_res} 
\end{align}

The bound of the first term in the last equality is easy. We start by writing
\eqh{
\frac\nu{\ep^2}\int_{\Omega_\res^\veps}{p'(\vr)\Grad b\cdot\Grad\frac\vr b}\,\dx\,=\frac\nu{\ep^2}\int_{\Omega_\res^\veps}{\sqrt{p'(\vr)\,\vr}\;\Grad\ln b\cdot\left(\sqrt{\frac{p'(\vr)}{\vr}}\,b\,\Grad\frac\vr b\right)}\,\dx.
}
Then, we observe that
\[
p'(\vr)\,\vr\,=\,\vr^\g+\vr^{-\kappa}\qquad\qquad\Longrightarrow\qquad\qquad \left[\sqrt{p'(\vr)\,\vr}\right]_\res\,\in\,L^\infty\big(\R_+;L^2\big)
\]
uniformly in $\veps>0$, owing to \eqref{ub:dens-res}. Indeed, one has $\left\|\left[\sqrt{p'(\vr)\,\vr}\right]_\res\right\|_{L^2}\,\leq \,\left\|\vr\right\|_{L^\g}^{\g/2}+\,
\left\|\vr^{-1}\right\|_{L^\kappa}^{\kappa/2}$.
Thus, using \eqref{ub:dens-res} quantitatively and arguing as in Case 1, we deduce that
\begin{align*}
\left|\frac\nu{\ep^2}\int_{\Omega_\res^\veps}{\sqrt{p'(\vr)\,\vr}\;\Grad\ln b\cdot\left(\sqrt{\frac{p'(\vr)}{\vr}}\,b\,\Grad\frac\vr b\right)}\,\dx\right|\,&\lesssim\,
\frac{C_\de}{\veps^2}\left\|\left[\sqrt{p'(\vr)\,\vr}\right]_\res\right\|^2_{L^2}\,+\,\delta\,\frac{\nu}{\veps^2}\intO{b^2\,\frac{p'(\vr)}{\vr}\,\left|\nabla\left(\frac{\vr}{b}\right)\right|^2} \\
&\lesssim\,C_\de\,+\,\delta\,\frac{\nu}{\veps^2}\intO{b^2\,\frac{p'(\vr)}{\vr}\,\left|\nabla\left(\frac{\vr}{b}\right)\right|^2}\,.
\end{align*}
Once again, after integration in time, the last term on the right can be absorbed in the left-hand side of \eqref{BD}, for $\de>0$ small enough.

It remains to deal with the last term appearing in \eqref{est:J4_res}.
We start by writing
\eqh{\frac\nu{\ep^2}\int_{\Omega_\res^\veps}{p'(b)\Grad b\cdot\Grad\frac\vr b}\,\dx
=\frac\nu{\ep^2}\int_{\Omega_\res^\veps}{\sqrt{\frac{\vr}{p'(\vr)}}\,\frac{p'(b)}{b}\Grad b\cdot\left(\sqrt{\frac{p'(\vr)}{\vr}}\,b\,\Grad\frac\vr b\right)}\,\dx\,.}
First note that
\[
\frac{\vr}{p'(\vr)}\,=\,\frac{\vr}{p_E'(\vr)+p'_c(\vr)}\,=\,\frac{\vr^{\k+2}}{\vr^{\g+\k}+1}\,.
\]
Thus, on the one hand we get
\[
\left|\frac\nu{\ep^2}\int_{\Omega_{\res,B}^\veps}{\sqrt{\frac{\vr}{p'(\vr)}}\,\frac{p'(b)}{b}\Grad b\cdot\left(\sqrt{\frac{p'(\vr)}{\vr}}\,b\,\Grad\frac\vr b\right)}\,\dx\right|\,\lesssim\,
\frac\nu{\ep^2}\,\Big(\mc L(\Omega^\veps_\res)\Big)^{1/2}\,\left\|\sqrt{\frac{p'(\vr)}{\vr}}\,b\,\Grad\frac\vr b\right\|_{L^2}\,.
\]
On the other hand, we also have
\[
\mathds{1}_{\Omega_{\res,UB}^\veps}\,\sqrt{\frac{\vr}{p'(\vr)}}\,\lesssim\,\mathds{1}_{\Omega_{\res,UB}^\veps}\,\vr_\veps^{(2-\g)/2}\,.
\]
Now, in view of \eqref{ub:dens-res}, for $ 2-\g\,\geq\,0\,$ we have that
\[
\left[\vr\right]_\res^{(2-\g)/2}\,\in\,L^\infty(\R_+;L^{q})\,,\quad q\,:=\,\frac{2\g}{2-\g}\,>\,2\,,\qquad\qquad \mbox{ with }\qquad
\left\|\left[\vr\right]_\res^{(2-\g)/2}\right\|_{L^q}\,=\,\left\|\left[\vr\right]_\res\right\|_{L^\g}^{\g/q}\,\lesssim\,\veps^{2/q}\,.
\]
If, on the other hand $ 2-\g\,<\,0$, then we have
\[
\left[\vr\right]_\res^{(2-\g)/2}\,\in\,L^\infty(\R_+;L^{q})\,,\quad q\,:=\,\frac{2\kappa}{\g-2}\,,\qquad\qquad \mbox{ with }\qquad
\left\|\left[\vr\right]_\res^{(2-\g)/2}\right\|_{L^q}\,=\,\left\|\left[\vr^{-1}\right]_\res\right\|_{L^\kappa}^{\kappa/q}\,\lesssim\,\veps^{2/q}\,.
\]
In the latter case, in order to have $q\geq2$, we need the condition
\[
\k\,\geq\,\g\,-\,2\,.
\]
In any case, and under the assumption that the previous requirement is satisfied in the case $\g>2$, we can thus bound
\[
\left|\frac\nu{\ep^2}\int_{\Omega_{\res,UB}^\veps}{\sqrt{\frac{\vr}{p'(\vr)}}\,\frac{p'(b)}{b}\Grad b\cdot\left(\sqrt{\frac{p'(\vr)}{\vr}}\,b\,\Grad\frac\vr b\right)}\,\dx\right|\,\lesssim\,
\frac{1}{\veps^2}\,\left\|\left[\vr\right]_\res^{(2-\g)/2}\right\|_{L^q}\,\left\|\sqrt{\frac{p'(\vr)}{\vr}}\,b\,\Grad\frac\vr b\right\|_{L^2}\,
\Big(\mc L(\Omega^\veps_\res)\Big)^{1/m}\,,
\]
where $1/q\,+\,1/m\,=\,1/2$.
In the end, using \eqref{ub:dens-res}, \eqref{ub:residual-set} and the definition of $m$, after an application of Young's inequality we arrive at
\[
\left|\frac\nu{\ep^2}\int_{\Omega_\res^\veps}{\sqrt{\frac{\vr}{p'(\vr)}}\,\frac{p'(b)}{b}\Grad b\cdot\left(\sqrt{\frac{p'(\vr)}{\vr}}\,b\,\Grad\frac\vr b\right)}\,\dx\right|\,\lesssim\,
C_\de\,+\,\de\,\frac{\nu}{\veps^2}\,\left\|\sqrt{\frac{p'(\vr)}{\vr}}\,b\,\Grad\frac\vr b\right\|^2_{L^2}\,.
\]

Finally, putting all those inequalities together, we see that
\[
\int^T_0\big|J_3\big|\,\dt\,\lesssim\,C_\de\,+\,3\,\de\,\frac{\nu}{\veps^2}\intTO{b^2\,\frac{p'(\vr)}{\vr}\,\left|\nabla\left(\frac{\vr}{b}\right)\right|^2}\,.
\]
Hence, for $\delta>0$ small enough, we can absorbe the last term of the previous inequality into the left-hand side of \eqref{BD}. The proof of the lemma is thus completed. \qed

Before moving on, let us observe that the constant in Lemma \ref{l:BD} depends on time. Therefore, all the bounds which we are going to derive
from the BD entropy estimate will be only \emph{local in time}, on any arbitrarily large but compact interval $[0,T]$.
This contrasts with the bounds coming from the classical energy inequality, which instead are global in time.

\subsection{Consequences of the BD entropy estimate} \label{ss:consequences}

Lemma \ref{l:BD} provides us with at least three additional pieces of information with respect to the classical energy inequality \eqref{est:energy}. Those additional bounds will be fundamental
in order to prove convergence in our setting.

First of all, estimate \eqref{BD_estimate} can be used to improve \eqref{ub:Du_zero}: indeed, as $\nabla\vu\,=\,\tD\vu\,+\,\mbb A\vu$, we have
\eq{
\big(\sqrt{\vr_\veps}\,\Grad \vu_\veps\big)_\veps\,&\subset\,L^2_{\rm loc}\big(\R_+;L^2(\Omega)\big)\,. \label{ub:Du} 
}

Next, combining \eqref{BD_estimate} with \eqref{est:energy}, we gather that
\begin{equation} \label{est:log}
\sqrt{\vr_\veps}\,\nabla\ln\left(\frac{\vr_\veps}{b}\right)\,=\,\frac{b}{\sqrt{\vr_\veps}}\,\nabla\left(\frac{\vr_\veps}{b}\right)\;\in\,L^\infty_{\rm loc}\big(\R_+;L^2(\Omega)\big)\,,
\end{equation}
with \emph{uniform} inclusion (of course, uniform with respect to $\veps\in\,]0,1]$).
At this point, an easy computation shows that
\[
\frac{b}{\sqrt{\vr_\veps}}\,\nabla\left(\frac{\vr_\veps}{b}\right)\,=\,2\,\nabla\sqrt{\vr_\veps}\,-\,\sqrt{\vr_\veps}\,\nabla\ln b\,=\,
2\,\nabla\sqrt{\vr_\veps}\,-\,\left(\sqrt{[\vr_\veps]_\ess}\,+\,\sqrt{[\vr_\veps]_\res}\right)\,\nabla\ln b\,.
\]
Observe that $\big([\vr_\veps]_\ess\big)_\veps$ is uniformly bounded in $L^\infty(\R_+\times\Omega)$, whereas $\left(\sqrt{[\vr_\veps]_\res}\right)_\veps$
is uniformly bounded in $L^\infty(\R_+;L^{2\g})\,\hookrightarrow\,L^\infty(\R_+;L^{2})$, in view of \eqref{ub:dens-res}.
Putting this information together with \eqref{est:log}, by Sobolev embeddings we finally deduce that
\begin{equation} \label{est:sqrt-rho}
\left(\nabla\sqrt{\vr_\veps}\right)_\veps\,\subset\,L^\infty_{\rm loc}\big(\R_+;L^2(\Omega)\big)\qquad\Longrightarrow\qquad \big(\vr_\veps\big)_\veps\,\subset\,L^\infty_{\rm loc}\big(\R_+;L^3(\Omega)\big)\,.
\end{equation}
The previous property provides higher integrability for $\vr_\veps$, \emph{globally} (\tsl{i.e.} in the whole $\Omega$, without having to distinguish between essential and residual parts) and
\emph{uniformly} with respect to $\veps\in\,]0,1]$.

Next, we consider the quantity
\begin{equation} \label{ub:Drho_BD}
\left(\frac{1}{\veps^2}\,{\frac{p'(\vr_\veps)}{\vr_\veps}}\,b^2\,\left|\Grad\left(\frac{\vr_\veps}{b}\right)\right|^2\right)_\veps\,\subset\,L^1_{\rm loc}\big(\R_+;L^1(\Omega)\big)\,.
\end{equation}
The term inside the parentheses can be written as
\begin{align*}
 {\frac{p'(\vr_\veps)}{\vr_\veps}}\,b^2\,\left|\Grad\left(\frac{\vr_\veps}{b}\right)\right|^2\,&=\,\lr{\vre^{\g-2}+\vre^{-\kappa-2}}\,b^2\,\left|\Grad\left(\frac{\vr_\veps}{b}\right)\right|^2\,\\
 &=\,
\left|b^{\g/2}\,\left(\frac{\vr_\veps}{b}\right)^{(\g-2)/2}\,\Grad\left(\frac{\vr_\veps}{b}\right)\right|^2\,
+\left|b^{-\kappa/2}\,\left(\frac{\vr_\veps}{b}\right)^{(-\kappa-2)/2}\,\Grad\left(\frac{\vr_\veps}{b}\right)\right|^2\,\\
&=\,\left|\frac{2}{\g}\,b^{\g/2}\,\nabla\left(\frac{\vr_\veps}{b}\right)^{\g/2}\right|^2\,+\left|-\frac{2}{\kappa}\,b^{-\kappa/2}\,\nabla\left(\frac{\vr_\veps}{b}\right)^{-\kappa/2}\right|^2\,.
\end{align*}
Thus, we get the following uniform embeddings:
\begin{equation} \label{est:Drho}
\left(\frac{1}{\veps}\,\nabla\left(\frac{\vr_\veps}{b}\right)^{\g/2}\right)_\veps\,\subset\,L^2_{\rm loc}\big(\R_+;L^2(\Omega)\big),\qquad\left(\frac{1}{\veps}\,\nabla\left(\frac{\vr_\veps}{b}\right)^{-\kappa/2}\right)_\veps\,\subset\,L^2_{\rm loc}\big(\R_+;L^2(\Omega)\big)\,.
\end{equation}
This is another fundamental piece of information, since it gives, roughly speaking, a uniform control on the gradient of $\vr_\veps/b$, with a quantitative bound
$O(\veps)$ in a suitable norm. However, in order to be able to fully exploit this information, some preparatory work is needed.

To begin with, we write
\[
\frac{1}{\veps}\,\nabla\left(\frac{\vr_\veps}{b}\right)^{\g/2}\,=\,\frac{1}{\veps}\nabla\left(\frac{\vr_\veps^{\g/2}\,-\,b^{\g/2}}{b^{\g/2}}\right)\,.
\]
Therefore, the uniform bound \eqref{est:Drho} implies that
\begin{equation} \label{est:conseq-1}
\left(\frac{1}{\veps}\,\left(\vr_\veps^{\g/2}\,-\,b^{\g/2}\right)\right)_\veps\,\subset\,L^2_{\rm loc}\big(\R_+;L^6(\Omega)\big)\,.
\end{equation}

Let us derive a couple of properties from \eqref{est:conseq-1}. We notice that we can write
\[
\vr_\veps^{\g/2}\,-\,b^{\g/2}\,=\,\frac{\vr_\veps^{\g}\,-\,b^{\g}}{\vr_\veps^{\g/2}\,+\,b^{\g/2}}\,.
\]

On the one hand, from the previous relation we immediately deduce that
\[
\left(\frac{1}{\veps}\,\mathds{1}_{\Omega^\veps_{\res,B}\cup\Omega^\veps_\ess}\,\left(\vr_\veps^{\g}\,-\,b^{\g}\right)\right)_\veps\,\subset\,L^2_{\rm loc}\big(\R_+;L^6(\Omega)\big)\,.
\]
We use a Taylor expansion of the function $f(s)\,=\,s^\g$ to get, for a suitable point $z_\veps\,=\,z_\veps(t,x)$ between $\vr_\veps(t,x)$ and $b(x)$, the following fact:
\[
 \vr_\veps^{\g}\,-\,b^{\g}\,=\,\g\,b^{\g-1}\,\big(\vr_\veps\,-\,b\big)\,+\,\g\,(\g-1)\,z_\veps^{\g-2}\,\big(\vr_\veps\,-\,b\big)^2\,\geq\,\g\,b^{\g-1}\,\big(\vr_\veps\,-\,b\big)\,.
\]
This inequality, together with \eqref{est:conseq-1} above, implies that
\begin{equation} \label{est:rho-B}
\left(\mathds{1}_{\Omega^\veps_{\res,B}\cup\Omega^\veps_\ess}\,\phi_\veps\right)_\veps\,\subset\,L^2_{\rm loc}\big(\R_+;L^6(\Omega)\big)\,.
\end{equation}

On the other hand, on the set $\Omega^\veps_{\res,UB}$ we have that $\vr_\veps\,\geq\,2\,b^*$. Hence, on that set we can write
\[
\frac{\vr_\veps^{\g}\,-\,b^{\g}}{\vr_\veps^{\g/2}\,+\,b^{\g/2}}\,\geq\,\left(\frac{\vr_\veps}{2}\right)^{\g}\,\frac{2^\g-1}{2\,\vr_\veps^{\g/2}}\,\geq\,C\,\vr_\veps^{\g/2}\,.
\]
Using \eqref{est:conseq-1} again, we gather that 
\begin{equation} \label{est:rho-UB}
\left(\frac{1}{\veps}\,\mathds{1}_{\Omega^\veps_{\res,UB}}\,\vr^{\g/2}_\veps\right)_\veps\,\subset\,L^2_{\rm loc}\big(\R_+;L^{6}(\Omega)\big)\qquad\quad\Longrightarrow\qquad\quad
\left(\frac{1}{\veps}\,\mathds{1}_{\Omega^\veps_{\res,UB}}\,\vr_\veps\right)_\veps\,\subset\,L^{\gamma}_{\rm loc}\big(\R_+;L^{3\g}(\Omega)\big)\,.
\end{equation}
The main point of the last computations is that \eqref{est:Drho} does not really provide a uniform control on the gradient of the functions $\phi_\veps$ in $L^2_T(L^2)$, whenever $\g\neq2$.
Such a control, which is true when $\g=2$, would give higher integrability of the $\phi_\veps$'s in $L^2_T(L^6)$.
Estimate \eqref{est:rho-B} shows that we are not too far from getting that property, 
 but this fact is true only when $\vr_\veps$ stays bounded.
Since we do not have $L^\infty$ bounds on $\vr_\veps$, the density functions may grow in some part of the domain $\Omega$ (which has to be small, recall \eqref{ub:dens-res} above).
Anyway, inequality \eqref{est:rho-UB} provides us with a useful control in that region.

\medbreak
Actually, in $\Omega^\veps_{\res,B}$ we can derive an even better control on the $\vr_\ep$'s than \eqref{est:rho-B}. Indeed, arguing in a similar way as for getting \eqref{est:rho-UB},
%In a very similar way,
from the second piece of information in \eqref{est:Drho} we deduce
\begin{equation} \label{est:rho-1}
\left(\frac{1}{\veps}\,\mathds{1}_{\Omega^\veps_{\res,B}}\,\vr_\veps^{-\k/2}\right)_\veps\,\subset\,L^2_{\rm loc}\big(\R_+;L^6(\Omega)\big)
\qquad\quad\Longrightarrow\qquad\quad
\left(\frac{1}{\veps}\,\mathds{1}_{\Omega^\veps_{\res,B}}\,\vr_\veps^{-1}\right)_\veps\,\subset\,L^\k_{\rm loc}\big(\R_+;L^{3\k}(\Omega)\big)\,.
\end{equation}

%\medbreak
We conclude this part by showing a third consequence of the BD entropy estimate \eqref{BD_estimate}, that is a control on the gradient functions $\Grad\phi_\veps$ in suitable Lebesgue norms. Since this is
a key property, we put it in the form of a proposition.

\begin{prop} \label{p:D-phi}
Let $\g>1$ and $\k>0$ such that $\kappa\geq \gamma-2$.
Let $\big(\vr_\veps,\vu_\veps\big)_\veps$ be a family of global in time finite energy weak solutions to system \eqref{main}, in the sense of Definition \ref{df:main},
related to initial data $\big(\vr_{0,\veps},\vu_{0,\veps}\big)_\veps$ such that \eqref{ub:data} holds. 

Then:
\begin{itemize}
 \item if $\g\geq2$, one has the uniform embedding $\;\big(\phi_\veps\big)_\veps\,\subset\,L^2_{\rm loc}\big(\R_+;H^1(\Omega)\big)$;
 \item when $\g<2$, one deduces instead $\;\big(\phi_\veps\big)_\veps\,\subset\,L^2_{\rm loc}\big(\R_+;W^{1,\g}(\Omega)\big)$.
\end{itemize}
\end{prop}

\pf
We start by observing that, as a consequence of \eqref{ub:Drho_BD}, we get in particular
\begin{equation} \label{est:D-phi_start}
\left({\vr_\veps}^{-(1-\g/2)}\,b\,\Grad\left(\frac{\phi_\veps}{b}\right)\right)_\veps\,\subset\,L^2_{\rm loc}\big(\R_+;L^2(\Omega)\big),\quad \left({\vr_\veps}^{-(1+\kappa/2)}\,b\,\Grad\left(\frac{\phi_\veps}{b}\right)\right)_\veps\,\subset\,L^2_{\rm loc}\big(\R_+;L^2(\Omega)\big)\,.
\end{equation}
Since on the set $\Omega^\veps_\ess\,\cup\,\Omega^\veps_{\res,B}$ we have $0\leq \vr_\veps\leq 2b^*$, we easily get
\[
\left|\mathds{1}_{\Omega^\veps_\ess\,\cup\,\Omega^\veps_{\res,B}}\;b\,\Grad\left(\frac{\phi_\veps}{b}\right)\right|\,\lesssim\,
\left|\mathds{1}_{\Omega^\veps_\ess\,\cup\,\Omega^\veps_{\res,B}}\;{\vr_\veps}^{-(1+\kappa/2)}\,b\,\Grad\left(\frac{\phi_\veps}{b}\right)\right|\,,
\]
which implies that, for all $T>0$ fixed, one has
\eq{\label{better_phi}
\left(\mathds{1}_{\Omega^\veps_\ess\,\cup\,\Omega^\veps_{\res,B}}\;b\,\Grad\left(\frac{\phi_\veps}{b}\right)\right)_\veps\,\subset\,L^2_T(L^2)\,.
}

On the subset $\Omega^\veps_{\res,UB}$, we need to consider two cases.

\medbreak
\noindent{\emph{Case 1.}}  For $1-\g/2<0$ we can employ the first part of \eqref{est:D-phi_start} and an argument similar to the one used above to deduce that
\[
\forall\;T\,>\,0\,,\qquad\qquad \left(\mathds{1}_{\Omega^\veps_{\res,UB}}\;b\,\Grad\left(\frac{\phi_\veps}{b}\right)\right)_\veps\,\subset\,L^2_T(L^2)\,.
\]

\noindent{\emph{Case 2.}}  For $1-\g/2\geq0$  instead, we can write
\begin{align*}
\left|\mathds{1}_{\Omega^\veps_{\res,UB}}\;b\,\Grad\left(\frac{\phi_\veps}{b}\right)\right|\,&=\,
\left|\mathds{1}_{\Omega^\veps_{\res,UB}}\;{\vr_\veps}^{-(1-\g/2)}\,b\,\Grad\left(\frac{\phi_\veps}{b}\right)\right|\;\big[\vr_\veps\big]_{\res,UB}^{1-\g/2}\,.
\end{align*}
Since, for any $T>0$ fixed, $\big(\big[\vr_\veps\big]_\res\big)_\veps$ is uniformly bounded in $L^\infty_T(L^\g)$, we have that
$\left(\big[\vr_\veps\big]_{\res,UB}^{1-\g/2}\right)_\veps$ is uniformly bounded in $L^\infty_T(L^{2\g/(2-\g)})$.
In turn, combining this information with \eqref{est:D-phi_start} implies that
\[
\forall\;T\,>\,0\,,\qquad\qquad \left(\mathds{1}_{\Omega^\veps_{\res,UB}}\;b\,\Grad\left(\frac{\phi_\veps}{b}\right)\right)_\veps\,\subset\,L^2_T(L^\g)\,.
\]
We have thus discovered that the family of $\Grad\big(\phi_\veps/b\big)$'s is uniformly bounded in the space $L^2_T(L^\g)$, for any time $T>0$.

In order to conclude, we simply write
\[
\Grad\left(\frac{\phi_\veps}{b}\right)\,=\,\frac{1}{b}\,\Grad\phi_\veps\,-\,\frac{1}{b^2}\,\phi_\veps\,\Grad b\,,
\]
and, when $\g\leq 2$, we bound $\big(\phi_\veps\big)_\veps$ uniformly in $L^\infty_T(L^\g)$, thanks to \eqref{ub:phi_ess} and \eqref{ub:phi_res}. When $\g>2$, instead,
we use \eqref{ub:phi_ess} again, together with the fact that, by definition of $\phi_\eps$, we have
\begin{align*}
\int_{\Omega_{\res,B}}\left|\phi_\eps\right|^2\,\dx\,&=\,\dfrac{1}{\veps^2}\int_{\Omega_{\res,B}}\left|\vr_\eps\,-\,b\right|^2\,\dx\;\lesssim\;\dfrac{1}{\veps^2}\,\mc L\big(\Omega_\res\big)\;\lesssim\;1 \\
\int_{\Omega_{\res,UB}}\left|\phi_\eps\right|^2\,\dx\,&=\,\dfrac{1}{\veps^2}\int_{\Omega_{\res,UB}}\left|\vr_\eps\,-\,b\right|^2\,\dx\;\lesssim\;
\dfrac{1}{\veps^2}\int_{\Omega_{\res,UB}}\vr_\eps^2\,\dx \\
&\lesssim\,\dfrac{1}{\veps^2}\,\left\|\big[\vr_\veps\big]_\res\right\|^2_{L^\g}\,\Big(\mc L\big(\Omega_\res\big)\Big)^{(\g-2)/\g}\;\lesssim\;1\,,
\end{align*}
where we have used also the bounds provided by \eqref{ub:dens-res}. We thus conclude that, when $\g>2$, we have $\big(\phi_\veps\big)_\eps\,\subset\,L^\infty_T(L^2)$ for any $T>0$ fixed.

In the end, the proposition is proved.
\qed

We also notice that, as a consequence of Proposition \ref{p:D-phi} and Sobolev embedding, we get that
\begin{equation} \label{ub:phi_higher}
\mbox{if }1<\g< 2\,,\quad \big(\phi_\veps\big)_\veps\,\subset\,L^2_{\rm loc}\big(\R_+;L^{3\g/(3-\g)}(\Omega)\big)\,,
\qquad\quad\mbox{and if }\g\geq2\,,\quad\big(\phi_\veps\big)_\veps\,\subset\,L^2_{\rm loc}\big(\R_+;L^{6}(\Omega)\big)\,.
\end{equation}

\section{The singular perturbation operator and the wave system} \label{s:singular}

In this section, we study in detail the singular part of the primitive equations \eqref{main}. 
To begin with, in Subsection \ref{ss:first-convergence} we establish first convergence properties for the sequence of solutions $\big(\vr_\veps,\vu_\veps\big)_\veps$
towards some targe profile $\big(\oline{\vr},\oline\vu\big)$. Then, in Subsection \ref{ss:constraints} we derive some constraints that the limit point $\big(\oline{\vr},\oline\vu\big)$
has to satisfy.
Finally, in Subsection \ref{ss:acoustic} we introduce the system of waves, which encodes the propagation of fast time oscillations.

\subsection{Preliminary convergence properties} \label{ss:first-convergence}

From the uniform bounds exhibited in Subsections \ref{ss:energy}, \ref{ss:BD} and \ref{ss:consequences}, we can derive first convergence properties for the family of solutions $\big(\vr_\veps,\vu_\veps\big)_\veps$
of our primitive system \eqref{main}. Of course, the convergence we are going to establish is only in weak topologies, therefore it will not be enough for deriving the limit system \eqref{eq:limit}.

Here below, it is convenient to work on time intervals $[0,T]$, for arbitrary large but fixed $T>0$.
Also, in the notation we imply that all the convergences are taken in the limit $\veps\ra0^+$.

\paragraph{The density functions.}
We start by considering the sequence of the density functions $\vr_\veps$. First of all, from \eqref{ub:dens-ess} and \eqref{ub:dens-res}, we see that
we can decompose
\begin{equation} \label{conv:rho}
\vr_\veps\,=\,b\,+\,\vr_\veps^{(1)}\,+\,\vr_\veps^{(2)}\,,\qquad \mbox{ with }\qquad
\vr_\veps^{(1)}\,\longrightarrow\,0 \ \mbox{ in }\ L^\infty_T(L^2)\,,\quad \vr_\veps^{(2)}\,\longrightarrow\,0\ \mbox{ in }\ L^\infty_T(L^\g)\,.
\end{equation}
Since  $\Omega$ is of finite measure, we can interpolate those convergence properties with \eqref{est:sqrt-rho} to deduce that
\[
\vr_\veps\,-\,b\,=\,r_\veps\,,\qquad\qquad \mbox{ with }\qquad\qquad r_\veps\,\longrightarrow\,0\qquad \mbox{ in }\quad L^\infty_T(L^p)\qquad \forall\;p<3\,.
\]

The uniform bounds of Section \ref{s:unif-bounds} allow us to find more quantitative convergence properties. Indeed, from \eqref{ub:phi_ess} and \eqref{est:rho-B}
it follows that there exists a function $\phi\in L^\infty\big(\R_+;L^2(\Omega)\big)\cap L^2_{\rm loc}\big(\R_+;L^6(\Omega)\big)$ such that,
up to the extraction of a suitable subsequence, one has
\begin{equation} \label{conv:phi}
\big[\phi_\veps\big]_\ess\,\stackrel{*}{\rightharpoonup}\,\phi\qquad\qquad\mbox{ in }\qquad L^\infty_T(L^2)\,\cap\,L^2_T(L^6)\,,
\end{equation}
where the symbol $\stackrel{*}{\rightharpoonup}$ stands for the weak-$*$ convergence in the respective functional space.
On the other hand, owing to \eqref{ub:phi_res}, we know that
\begin{equation} \label{conv:phi-res}
\big[\phi_\veps\big]_\res\,\longrightarrow\,0\qquad\qquad \mbox{ in }\qquad L^\infty_T(L^p)\,,\qquad \forall\,p\in[1,\min\{\g,2\}[\,.
\end{equation}
Of course, if $\g\geq2$, we have the previous strong convergence only for all $1\leq p<2$; however, interpolating with the bounds of \eqref{ub:phi_higher}, in that case we get
\[
\big[\phi_\veps\big]_\res\,\longrightarrow\,0\qquad\qquad \mbox{ in }\qquad L^2_T(L^p)\,,\qquad \forall\,p\in[1,6[\,.
\]

\paragraph{The velocity fields.}
As it is apparent from equations \eqref{main}, any information on the velocity fields $\vu_\veps$ and their gradients is lost in regions close to vacuum.
This is one of the main difficulties arising in the analysis of system \eqref{main}. On the other hand,
at least at a first sight, it is not so clear which is the right quantity to look at; for instance, keep in mind inequalities \eqref{ub:kinetic}, \eqref{ub:Du}. 
Here, we decide to work with the momentum
\[
\vc V_\veps\,:=\,\vr_\veps\,\vu_\veps\,.
\]
However, the first step is to get some uniform bounds on the velocity fields $\vu_\veps$. This is the goal of the next proposition.
\begin{prop} \label{p:u}
Let $\big(\vr_\veps,\vu_\veps\big)_\veps$ be a family of global in time finite energy weak solutions to system \eqref{main}, in the sense of Definition \ref{df:main},
related to initial data $\big(\vr_{0,\veps},\vu_{0,\veps}\big)_\veps$ such that \eqref{ub:data} holds.

Then we have the uniform estimates
\[
\big(\vu_\eps\big)_\eps\,\subset\,L^\infty\big(\R_+;L^{p_1}(\Omega)\big)\qquad\mbox{ and }\qquad
\big(\Grad\vu_\eps\big)_\eps\,\subset\,L^2_{\rm loc}\big(\R_+;L^{p_1}(\Omega)\big)\,,\qquad\qquad p_1\,:=\,\frac{2\k}{\k+1}\,.
\]
In particular, we also have that 
\[
\big(\vu_\eps\big)_\eps\,\subset\,L^2_{\rm loc}\big(\R_+;L^{p_2}(\Omega)\big)\,,\qquad\qquad\qquad\mbox{ where }\qquad p_2\,:=\,\frac{6\k}{\k+3}\,.
\]

\end{prop}
\begin{rmk}\label{rem:p1p2}
Note that, due to our assumption $\kappa>3$, one has $p_1>3/2$ and $p_2>3$ in the above proposition.
\end{rmk}

\pf
We start by writing $\vue\,=\,\vre^{-1/2}\,\sqrt{\vre}\,\vue$, from which, by use of \eqref{ub:kinetic} and \eqref{ub:dens-res}, we immediately deduce
that $\big(\vue\big)_\eps$ is bounded in $L^\infty_T(L^{p_1})$.

For the gradient terms, the argument is analogous, since we can write
$\Grad\vue\,=\,\vre^{-1/2}\,\sqrt{\vre}\,\Grad\vue$ and use, this time, \eqref{ub:Du} and \eqref{ub:dens-res}.
This implies the claimed bound on the sequence $\big(\Grad\vue\big)_\veps$.

Finally, the last uniform bound for the family of $\vue$'s follows from the previous property and Sobolev embeddings.
\qed

We are now ready to derive some important estimates for the velocity fields $\vc V_\veps$.
\begin{prop} \label{p:momentum}
Let $\big(\vr_\veps,\vu_\veps\big)_\veps$ be a family of global in time finite energy weak solutions to system \eqref{main}, in the sense of Definition \ref{df:main},
related to initial data $\big(\vr_{0,\veps},\vc m_{0,\veps}\big)_\veps$ such that \eqref{ub:data} holds. For all $\veps\in\,]0,1]$, define $\vc V_\veps\,:=\,\vr_\veps\,\vu_\veps$.

Then, the following facts hold:
\begin{enumerate}[(i)]
 \item the sequence $\big(\vc V_\eps\big)_\eps$ is uniformly bounded in the space $L^\infty_{\rm loc}\big(\R_+;L^{3/2}(\Omega)\big)\,\cap\,L^2_{\rm loc}\big(\R_+;W^{1,p_3}(\Omega)\big)$,
where $p_3:=6\kappa/(5\kappa+3)$;
 \item there exist sequences $\big(\mc V_\eps\big)_\eps$ and $\big(\vc W_\eps\big)_\eps$ of vector fields such that 
\[ %\begin{equation} \label{eq:mom-decomp_ep}
\forall\,\veps\in\,]0,1]\,,\qquad\vc V_{\eps}\,=\,\mc V_{\eps}\,+\,\eps\,\vc W_{\eps,M}\,,
\] %\end{equation}
with the uniform embedding properties
\[
\big(\mc V_{\eps}\big)_{\eps}\,\subset\,L^\infty_{\rm loc}\big(\R_+;L^2(\Omega)\big)\qquad \mbox{ and }\qquad \big(\vc W_{\eps}\big)_{\eps}\,\subset\,L^2_{\rm loc}\big(\R_+;L^{3/2}(\Omega)\big)\,;
\]
\item after setting $\wtilde{\mc V}_\eps\,:=\,b\,\vu_\ep$, we can also write
\[ %\begin{equation} \label{eq:mom-decomp2_ep}
\vc V_{\eps}\,=\,\wtilde{\mc V}_{\eps}\,+\,\eps\,\wtilde{\vc W}_{\eps}\,, %\qquad\qquad \wtilde{\mc V}_\eps\,:=\,b\,\vu_\ep\,,
\] %\end{equation}
with the uniform embedding properties
\[
\big(\wtilde{\mc V}_{\eps}\big)_{\eps}\,\subset\,L^2_{\rm loc}\big(\R_+;W^{1,p_1}(\Omega)\big)\qquad \mbox{ and }\qquad
\big(\wtilde{\vc W}_{\eps}\big)_\eps\,\subset\,L^2_{\rm loc}\big(\R_+;L^{1}(\Omega)\big)\,. \]
\end{enumerate}
\end{prop}

\pf
Let $T>0$ be arbitrary, but fixed throughout the following computations.

The proof of the $L^\infty_T(L^{3/2})$ bound of item (i) is easy to get: it is enough to write $\vc V_\veps\,=\,\sqrt{\vr_\veps}\;\sqrt{\vr_\veps}\,\vu_\veps$ and use the uniform bounds given in
\eqref{ub:kinetic} and \eqref{est:sqrt-rho}.
Next, let us focus on the bounds for the gradient. For any $1\leq j\leq 3$, we compute
\[
\d_j\vc V_\eps\,=\,\sqrt{\vr_\eps}\,\sqrt{\vr_\eps}\,\d_j\vu_\eps\,+\,\vu_\eps\,\d_j\vr_\eps\,=\,
\sqrt{\vr_\eps}\,\sqrt{\vr_\eps}\,\d_j\vu_\eps\,+\,2\,\sqrt{\vr_\veps}\,\vu_\eps\,\d_j\sqrt{\vr_\eps}\,:=\,A_\eps\,+\,B_\eps\,.
\]
Repeating the previous argument, using this time \eqref{ub:Du} and \eqref{est:sqrt-rho}, it is easy to see that $\big(A_\eps\big)_\eps\,\subset\,L^2_T(L^{3/2})$. As for $B_\eps$, we start by observing that,
owing to Proposition \ref{p:u}, we have
$\big(\vu_\eps\big)_\eps\,\subset\,L^2_T(L^p)$ for $p\leq p_2:= \frac{6\kappa}{\kappa+3}$.
On the other hand, we have $\big(\Grad\sqrt{\vr_\eps}\big)_\eps\,\subset\,L^\infty_T(L^2)$ and $\big(\sqrt{\vr_\eps}\big)_\veps\,\subset\,L^\infty_T(L^6)$. Hence, we get
$\big(B_\eps\big)_\eps\,\subset\,L^2_T(L^{p_3})$ where $p_3=\frac{6\kappa}{5\kappa+3}$. Notice that $p_3>1$ if and only if $\kappa>3$:
this is precisely the place where the strongest assumption on $\kappa$ appears.
Item (i) is thus proven.

For showing the decomposition in item (ii), we notice that
\[
 \vc V_\veps\,=\,\sqrt{b}\,\sqrt{\vr_\eps}\,\vu_\eps\,+\,\eps\,\frac{\sqrt{\vr_\eps}\,-\,\sqrt{b}}{\eps}\,\sqrt{\vr_\eps}\,\vu_\eps\,.
\]
Thus, if we define $\mc V_\eps\,:=\,\sqrt{b}\,\sqrt{\vr_\eps}\,\vu_\eps$, from \eqref{ub:kinetic} we immediately infer that $\big(\mc V_{\eps}\big)_\eps\,\subset\,L^\infty_T(L^2)$. Next, we define
\[
\vc W_\eps\,:=\,\frac{\sqrt{\vr_\eps}\,-\,\sqrt{b}}{\eps}\,\sqrt{\vr_\eps}\,\vu_\eps\,.
\]
At this point, we use that $\left|\sqrt{\vr_\eps}-\sqrt{b}\right|\,\leq\,b_*^{-1/2}\,\left|\vr_\eps-b\right|$. Hence, on the one hand, by using \eqref{est:rho-B}, we see that
\[
\left(\mathds{1}_{\Omega^\veps_{\res,B}\cup\Omega^\veps_\ess}\,\frac{\sqrt{\vr_\eps}\,-\,\sqrt{b}}{\eps}\right)_\veps\,\subset\,L^2_T(L^6)\,;
\]
on the other hand, since
\[
0\,\leq\,\mathds{1}_{\Omega^\veps_{\res,UB}}\,\left(\sqrt{\vr_\eps}-\sqrt{b}\right)\,\leq\,\mathds{1}_{\Omega^\veps_{\res,UB}}\,\sqrt{\vr_\veps}\,\lesssim\,
\mathds{1}_{\Omega^\veps_{\res,UB}}\,\vr_\veps^{\g/2}\,,
\]
in view of \eqref{est:rho-UB} we gather that also $\left(\mathds{1}_{\Omega^\veps_{\res,UB}}\,\frac{\sqrt{\vr_\eps}\,-\,\sqrt{b}}{\eps}\right)_\veps$
is uniformly bounded in $L^2_T(L^6)$. Therefore, we get that
\[
 \left(\frac{\sqrt{\vr_\eps}\,-\,\sqrt{b}}{\eps}\right)_\veps\,\subset\,L^2_T(L^6)\qquad\quad \Longrightarrow\qquad\quad
 \big(\vc W_\veps\big)_\veps\,\subset\,L^2_T(L^{3/2})\,.
\]

The proof of item (iii) is similar. This time, we write $\vre-b\,=\,\big(\sqrt\vre-\sqrt b\big)\big(\sqrt\vre+\sqrt b\big)$ and get
\[
\vc V_{\veps}\:=\,b\,\vu_\veps\,+\,\big(\sqrt\vre-\sqrt b\big)\,\sqrt{\vre}\,\vue\,+\,\big(\sqrt\vre-\sqrt b\big)\,\sqrt b\,\vue\,.
\]
If we set $\wtilde{\mc V}_\eps\,:=\,b\,\vue$, Proposition \ref{p:u} ensures us that the claimed bounds for this quantity are satisfied. Next, we claim that
the sequence of
\[
\wtilde{\vc W}_\eps\,:=\,\frac{\sqrt\vre-\sqrt b}{\ep}\,\sqrt{\vre}\,\vue\,+\,\frac{\sqrt\vre-\sqrt b}{\ep}\,\sqrt b\,\vue
\]
is uniformly bounded in $L^2_T(L^1)$, for all $T>0$ fixed. For this, owing to \eqref{ub:kinetic}, the last part of Proposition \ref{p:u} and Remark \ref{rem:p1p2},
it is enough to show that
\begin{equation} \label{est:to-prove}
\sqrt{\vre}\,-\,\sqrt{b}\quad\mbox{ is of order }\;\; O(\veps)\qquad\qquad \mbox{ in }\qquad L^\infty_T(L^2)\,.
\end{equation}
As a matter of fact, since $\left|\sqrt{\vr_\eps}-\sqrt{b}\right|\,\leq\,b_*^{-1/2}\,\left|\vr_\eps-b\right|$, from \eqref{ub:dens-ess} and \eqref{ub:dens-res} we immediately get
that $\left(\left(\mathds{1}_{\Omega^\eps_\ess}\,+\,\mathds{1}_{\Omega^\eps_{\res,B}}\right)\frac{\sqrt{\vr_\eps}-\sqrt{b}}{\ep}\right)_\veps$
is uniformly bounded in $L^\infty_T(L^2)$. Finally, as done above, we have
\[
0\,\leq\,\mathds{1}_{\Omega^\veps_{\res,UB}}\,\left(\sqrt{\vr_\eps}-\sqrt{b}\right)\,\leq\,\mathds{1}_{\Omega^\veps_{\res,UB}}\,\sqrt{\vr_\veps}\,,
\]
which, in view of \eqref{ub:dens-res} again, implies that
\begin{align*}
\left\|\mathds{1}_{\Omega^\veps_{\res,UB}}\,\left(\sqrt{\vr_\eps}-\sqrt{b}\right)\right\|_{L^2}\,&\lesssim\,\left\|\mathds{1}_{\Omega^\veps_{\res,UB}}\,\sqrt{\vr_\veps}\right\|_{L^{2\g}}\,
\Big(\mc L(\Omega^\eps_\res)\Big)^{1/2-1/2\g}\,\lesssim\,\veps\,. 
\end{align*}
%where we have used \eqref{ub:dens-res} again.
In the end, we have shown \eqref{est:to-prove}, which in turn implies the sought bound for the vector fields $\wtilde{\vc W}_\eps$.
\qed

\bigbreak
From the previous proposition, we immediately deduce the following corollary. The proof is rather straightforward, hence omitted.
\begin{coroll} \label{c:momentum}
Under the assumptions of Proposition \ref{p:momentum}, there exists a vector field $\vc V$, belonging to $L^\infty_{\rm loc}\big(\R_+;L^{3/2}(\Omega)\big)\,\cap\,L^2_{\rm loc}\big(\R_+;W^{1,p_3}(\Omega)\big)$
such that, up to the extraction of a suitable subsequence, one has $\vc V_\veps\,\stackrel{*}{\rightharpoonup}\,\vc V$ in that space.

In addition, $\vc V$ also belongs to $L^\infty_{\rm loc}\big(\R_+;L^2(\Omega)\big)\,\cap\,L^2_{\rm loc}\big(\R_+;W^{1,p_1}(\Omega)\big)$, and, up to further extractions,
for all $T>0$ one has that $\mc V_\veps\,\stackrel{*}{\rightharpoonup}\,\vc V$ in $L^\infty_T(L^2)$ and $\wtilde{\mc V}_\veps\,\rightharpoonup\,\vc V$ in $L^2_T(W^{1,p_1})$. %for $q={\frac{6\kappa}{\kappa+3}}$.
\end{coroll}

To conclude this part, we define the target velocity field $\vc U$ as
\begin{equation} \label{def:u}
\vc U\,:=\,\frac{1}{b}\,\vc V\,,
\end{equation}
where $\vc V$ is the vector field identified in Corollary \ref{c:momentum}.

We notice that $\vc U$ is (up to a further extraction) the  weak-limit point of the sequence $\big(\vu_\veps\big)_\veps$ in the functional spaces identified in Proposition \ref{p:u}.
We point out also that $\vc U$ belongs to the same functional spaces to which $\vc V$ belongs.

\subsection{Constraints on the limit} \label{ss:constraints}

In the previous subsection, we have identified the limit points $b$, $\phi$ and $\vc V$ of the families $\big(\vr_\veps\big)_\veps$, $\big(\phi_\veps\big)_\veps$ and $\big(\vc V_\veps\big)_\veps$,
respectively.
Our next goal is to find some properties those limit points have to satisfy. We point out that these conditions do not fully characterise the limit dynamics, which will be deduced
by passing to the limit in the momentum equation.

\begin{prop} \label{p:constraints}
Let $\big(\vr_\veps,\vu_\veps\big)_\veps$ be a family of global in time finite energy weak solutions to system \eqref{main}, in the sense of Definition \ref{df:main},
related to initial data $\big(\vr_{0,\veps},\vu_{0,\veps}\big)_\veps$ such that \eqref{ub:data} holds. For all $\veps\in\,]0,1]$, define $\vc V_\veps\,:=\,\vr_\veps\,\vu_\veps$.
Let $\phi$ be the scalar function identified in \eqref{conv:phi}, $\vc V$ the vector field identified in Corollary \ref{c:momentum} and $\vc U$ the vector field defined in \eqref{def:u}.

Then $\vc U$ has to satisfy the \emph{anelastic constraint}
\[
\div\big(b\,\vc U\big)\,=\,0\qquad\qquad\mbox{ in }\qquad \mc D'\big(\R_+\times\Omega\big)\,.
\]
On the other hand, $\phi\,=\,\phi(b)$ is determined as a function of $b$ only.
\end{prop}

Before proving the previous proposition, let us state a simple lemma. It will be useful to understand how to deal with the singularity of the pressure and gravitational terms, and
to put in evidence the right singularity in the momentum equation.
\begin{lemma} \label{l:pressure}
Let the assumptions of Proposition \ref{p:constraints} be in force. Then, after defining
\[
%\begin{equation} \label{def:PI}
\Pi(\vr_\veps;b)\,:=\,p(\vr_\eps)\,-\,p(b)\,-\,p'(b)\,\big(\vr_\veps\,-\,b\big)\,,
%\end{equation}
\]
we have the equality
\[
\frac{1}{\veps^2}\,\Big(\nabla p(\vr_\eps)\,-\,\vr_\ep\nabla G\Big)\,=\,\frac{1}{\veps}\,b\,\nabla\big(H''(b)\,\phi_\ep\big)\,+\,\frac{1}{\veps^2}\nabla\Pi(\vr_\ep;b)\,.
\]
In addition, one has that $\big(\frac{1}{\veps^2}\,\Pi(\vr_\veps;b)\big)_\veps$ is uniformly bounded in $L^\infty\big(\R_+;L^1(\Omega)\big)$.
\end{lemma}

\pf 
The claimed identity follows from a simple algebraic computation:
\begin{align*}
\frac{1}{\veps^2}\,\Big(\nabla p(\vr_\eps)\,-\,\vr_\ep\nabla G\Big)\,&=\,\frac{1}{\veps^2}\,\Big(\nabla\big( p(\vr_\eps)\,-\,p(b)\big)\,-\,\big(\vr_\ep\,-\,b\big)\,\nabla G\Big) \\
&=\,\frac{1}{\veps^2}\nabla\Pi(\vr_\veps;b)\,+\,\frac{1}{\veps^2}\,\bigg(\nabla\Big(p'(b)\,\big(\vr_\veps\,-\,b\big)\Big)\,-\,\big(\vr_\ep\,-\,b\big)\,\nabla G\bigg) \\
&=\,\frac{1}{\veps}\,b\,\nabla\big(H''(b)\,\phi_\ep\big)\,+\,\frac{1}{\veps^2}\nabla\Pi(\vr_\ep;b)\,,
\end{align*}
where, in the last step, we have used the equality $\nabla G\,=\,H''(b)\,\nabla b$, which holds owing to \eqref{eq:b}. %, we have
%\[
%\nabla G\,=\,H''(b)\,\nabla b\,.
%\]

It remains us to show the uniform bounds for the family $\big(\frac{1}{\veps^2}\,\Pi(\vr_\veps;b)\big)_\veps$. The argument is similar to the proof of Lemma 4.1 in \cite{F_2019}; however,
since that paper did not deal with non-constant limit density profiles $b$ nor with the presence of a cold component $p_c$ of the pressure, we report here most of the details.

First of all, using Taylor formula at the second order for the pressure function $p$ together with \eqref{ub:phi_ess}, we easily see that
\[
\frac{1}{\veps^2}\left[\Pi(\vr_\ep;b)\right]_\ess\,\approx\,\phi_\ep^2\qquad\qquad\Longrightarrow\qquad\qquad
\left(\frac{1}{\veps^2}\left[\Pi(\vr_\ep;b)\right]_\ess\right)_\eps\,\subset\,L^\infty(\R_+;L^1)\,.
\]
Next, let us denote by $\Pi_E(\vr_\ep;b)$ and $\Pi_c(\vr_\ep;b)$ the functions defined as $\Pi(\vr_\ep;b)$, but using respectively  $p_E$ and $p_c$ instead of the full pressure function $p$.
Then, a Taylor expansion again and \eqref{ub:residual-set} imply that
\[
\frac{1}{\veps^2}\left(\big[\Pi_E(\vr_\ep;b)\big]_{\res,B}\,+\,\big[\Pi_c(\vr_\ep;b)\big]_{\res,UB}\right)\,\subset\,L^\infty(\R_+;L^1)\,.
\]
Finally, we notice that
\[
\left|\big[\Pi_E(\vr_\ep;b)\big]_{\res,UB}\right|\,\lesssim\,\big[\vr_\ep\big]^\g_{\res}\qquad\qquad\mbox{ and }\qquad\qquad
\left|\big[\Pi_c(\vr_\ep;b)\big]_{\res,B}\right|\,\lesssim\,\big[\vr_\ep\big]^{-\k}_{\res}\,,
\]
for which one can deduce the sought bounds by using the controls in \eqref{ub:dens-res}.

In the end, the lemma is completely proved. \qed

We can now turn our attention to the proof of Proposition \ref{p:constraints}.

\medbreak
\noindent {\it Proof of Proposition \ref{p:constraints}.}
We start by considering the weak form of the momentum equation: given a test function $\vphi\in\mc D\big(\R_+\times\Omega\big)$, with $\Supp\vphi\subset[0,T[\,\times\Omega$ for some positive
time $T>0$, we have
\[
-\intTO{\vr_\veps\,\d_t\vphi}\,-\,\intTO{\vr_\veps\,\vu_\veps\cdot\nabla\vphi}\,=\,\intO{\vr_{0,\veps}\,\vphi(0)}\,.
\]
By assumption on the initial data, we know that $\vr_{0,\veps}\,-\,b\,\longrightarrow\,0$ in $L^2\cap L^\infty$, so it is easy to pass to the limit in the right-hand side of the previous relation.
Moreover, in view of \eqref{conv:rho} and Corollary \ref{c:momentum}, we know that $\vr_\veps\,\longrightarrow\,b$ and $\vr_\veps\,\vu_\veps\,=\,\vc V_\veps\,\longrightarrow\,\vc V$
in $\mc D'([0,T[\,\times\Omega)$. Thus, passing to the limit in the previous equality yields the constraint
\[
\intTO{\vc V\cdot\nabla\vphi}\,=\,0\qquad\forall\,\vphi\in\mc D\qquad\qquad\Longrightarrow\qquad\qquad \div \vc V\,=\,0\qquad \mbox{ in }\quad \mc D'\,.
\]
Now, using the definition of $\vc U$ given in \eqref{def:u} immediately gives the anelastic constraint.

\medbreak
Let us turn our attention to the momentum equation. As appears in the statement, the momentum equation does not give any relevant information on the limit points: we will discover
that $\phi$ is a quantity which play no role in the limit dynamics.

Indeed, we can test the momentum equation by $\veps\,\psi$, where $\psi\in\mc D(\R_+\times\Omega;\R^3)$ is a test function with support (say) included in $[0,T[\,\times\Omega$,
for some $T>0$. By the uniform bounds established in Section \ref{s:unif-bounds}, and using the identity in Lemma \ref{l:pressure} above, it is possible to see
that all the terms in the momentum equation tend to $0$, except
for the term $b\,\nabla\big(H''(b)\,\phi_\ep\big)$. Thus, passing to the limit for $\ep\ra0^+$, we find that
\[
\intTO{H''(b)\,\phi\,\div\big(b\,\psi\big)}\,=\,0\qquad\qquad\forall\,\psi\in\mc D\,.
\]
Since $b$ never vanishes, this implies that $\nabla\big(H''(b)\,\phi\big)\,=\,0$, and then $H''\big(b(x)\big)\,\phi(t,x)\,=\,c(t)$ a.e. on $\R_+\times\Omega$,
for a suitable function $c(t)$ only depending on the time variable.
We claim that $c(t)\,\equiv\,c(0)$ is in fact constant, and does not depend on time neither.
Indeed, let us denote by $\lan f\ran$ the space average over $\Omega$ of a function $f=f(t,x)$, namely
\[
\lan f\ran\,=\,\lan f\ran(t)\,:=\,\frac{1}{\mc L(\Omega)}\,\intO{f(t,x)}\,.
\]
Using the definition of $\phi_\veps$, we can recast the continuity equation as
\begin{equation} \label{eq:cont_wave}
\d_t\phi_\veps\,+\,\frac{1}{\veps}\,\div \vc V_\veps\,=\,0\,. 
\end{equation}
Taking the mean value with respect to the space variable, we discover that $\d_t\lan\phi\ran\,=\,0$,
thus the mean value of $\phi$ is preserved in time: $\lan\phi\ran(t)\,=\,\lan\phi_0\ran$ for almost all times $t\geq0$, where $\phi_0$ is the weak limit point of the initial data $\big(\phi_{0,\veps}\big)_\ep$
specified in \eqref{data_limit}.
On the other hand, coming back to the equality $H''(b)\,\phi\,=\,c(t)$ and computing the space average, we have that $\lan\phi\ran(t)\,=\,\lan1/H''(b)\ran\,c(t)$,
which in turn implies that $c(t)\equiv c_0$ has to be a constant function.

The proposition is now proved. \qed

%\bigbreak
As already pointed out, the previous proposition does not allow us to identify the limit dynamics yet. The main problem consists in passing to the limit in the momentum equation,
which reveals a not so easy task, owing to the non-linear terms appearing in it. In order to succeed, we first need to understand the propagation of fast oscillations in time: this 
is the goal of the next subsection.

\subsection{Acoustic equation} \label{ss:acoustic}

Lemma \ref{l:pressure} and Proposition \ref{p:constraints} allow us to identify the singular part of the equations of motion. 
Because of the ill-preparation of the initial data, this singular part is responsible for fast oscillations in time
of the solution, the so called \emph{acoustic waves}, which may eventually prevent the convergence of the non-linear terms. 
In order to study those oscillations and be able to pass to the limit in the equations, we reformulate our system  \eqref{main} as a wave system.

Let us recall that we have denoted
\begin{equation*}
\phi_\ep\,=\,\frac{\vr_\ep\,-\,b}{\ep}\qquad\qquad\mbox{ and }\qquad\qquad \vc V_\veps\,=\,\vr_\veps\,\vu_\veps\,.
\end{equation*}
As already remarked in the proof of Proposition \ref{p:constraints}, in terms of those quantities the continuity equation can be rewritten as in \eqref{eq:cont_wave}.
Concerning the momentum equation, instead, we can take advantage of the identity of Lemma \ref{l:pressure} to combine the pressure and gravitation terms together.

In the end, we discover that system \eqref{main} can be recasted in the following form:
\begin{equation} \label{eq:wave}
\left\{ \begin{array}{l}
 \ep\,\pt\phi_{\ep}\,+\,\Div \vc V_{\ep}\,=\,0 \\[1ex]
\ep\,\pt \vc V_{\ep}\,+\,b\,\Grad\big( H''(b)\, \phi_{\ep}\big)\,=\,\ep\,\vc f_{\ep}\,,
        \end{array}
\right.
\end{equation}
where we have defined
\begin{equation} \label{def:f_e}
\vc f_\ep\,=\,\nu\,\Div (\vr_\ep\,\tD\vu_\ep)\,-\,\Div (\vr_\eps \vu_\eps \otimes \vu_\eps)\,-\,\frac{1}{\veps^2}\Grad\Pi(\vr_\veps;b)\,.
\end{equation}

By the uniform bounds \eqref{ub:kinetic}, we get that the sequence $\big(\vr_\eps \vu_\eps \otimes \vu_\eps\big)_\veps$ is  uniformly bounded in $L^\infty_T(L^1)$, for all times $T>0$.
On the other hand, by arguing as in the proof of Proposition \ref{p:momentum}, it is easy to see that $\big(\vr_\veps\,\tD\vu_\veps\big)_\veps$ is uniformly bounded in $L^2_T(L^{3/2})$.
Finally, the term $\Pi(\vr_\veps;b)$ has been dealt with in Lemma \ref{l:pressure}.
Therefore, owing to the Sobolev embedding $H^s(\Omega)\subset L^\infty(\Omega)$ for any $s>3/2$, we get that
\begin{equation} \label{unif-b:f}
 \forall\,T>0\,,\quad \forall\,s>5/2\,,\qquad\qquad \big(\vc f_\eps\big)_\eps\,\subset\,L^2_T(H^{-s})\,.
\end{equation}

%\bigskip
For later use, it is convenient to introduce a regularised version of the wave system \eqref{eq:wave}. For this, we employ the low frequency cut-off operator $S_M$, with $M\in\N$,
introduced in relation \eqref{eq:S_j} in the Appendix.

Since $S_M$ commutes with all partial derivatives, applying operator $S_M$ to \eqref{eq:wave} yields
\begin{equation} \label{eq:wave_M}
\left\{ \begin{array}{l}
 \ep\,\pt\phi_{\ep,M}\,+\,\Div \vc V_{\ep,M}\,=\,0 \\[1ex]
\ep\,\pt \vc V_{\ep,M}\,+\,b\,\Grad\big( H''(b)\, \phi_{\ep,M}\big)\,=\,\ep\, \vc f_{\ep,M}\,+\,\vc h_{\ep,M}\,,
        \end{array}
\right.
\end{equation}
where we have denoted $\phi_{\veps,M}\,=\,S_M(\phi_\veps)$, $\vc V_{\veps,M}\,=\,S_M(V_\veps)$ and $\vc f_{\veps,M}\,=\,S_M(\vc f_\veps)$, and we have set
\[
 \vc h_{\ep,M}\,:=\,\big[b,S_M\big]\Grad\big( H''(b) \phi_{\ep}\big)\,+\,b\,\Grad\big( \big[H''(b),S_M\big] \phi_{\ep}\big)\,.
\]

Thanks to the uniform bounds for the sequence $\big(\phi_\eps\big)_\veps$ in $L^2_{\rm loc}\big(\R_+;W^{1,\g}(\Omega)\big)$ if $\g<2$, in $L^2_{\rm loc}\big(\R_+;H^1(\Omega)\big)$ when $\g\geq2$ (keep in mind
Proposition \ref{p:D-phi} above), standard commutator estimates (see \tsl{e.g.} Lemma 2.97 of \cite{B-C-D}) imply that
\begin{equation} \label{est:h}
\forall\;T\,>\,0\,,\qquad\qquad \left\|\vc h_{\ep,M}\right\|_{L^2_T(L^\g)}\,+\,\left\|\curl \vc h_{\ep,M}\right\|_{L^2_T(L^\g)}\,\lesssim\,2^{-M}\,,
\end{equation}
where we agree that the Lebesgue exponent $\g$ is changed into $2$ whenever $\g\geq2$.

We explicitly point out that, in the above estimate \eqref{est:h}, the multiplicative constant is uniform with respect to \emph{both} $M\in\N$ and $\eps\in\,]0,1]$, but it may depend on the fixed time $T>0$.
Notice that, for the uniform bound on $\curl \vc h_{\veps,M}$, the gradient structure of the commutator terms play a key role. We remark also that the commutator term $\vc h_{\veps,M}$
is much better controlled than the corresponding term in \cite{Masm}, thanks to the uniform bounds provided by the BD entropy estimates.

In addition, as an immediate consequence of \eqref{unif-b:f}, we get that
\begin{equation} \label{est:f}
\forall\;T\,>\,0\,,\quad \forall\;s\,\geq\,0\,,\qquad\qquad \left\|\vc f_{\ep,M}\right\|_{L^2_T(H^s)}\,\leq\,C(s,M)\,,
\end{equation}
where the constants $C(s,M)>0$ depend only on the quantities inside the parentheses. Notice that these constants blow up when $M\ra+\infty$, but they have finite value for any $M\in\N$ fixed.

We conclude this part by showing uniform bounds on the two sequences $\big(\phi_{\veps,M}\big)_{\veps,M}$ and $\big(V_{\veps,M}\big)_{\veps,M}$. As a matter of fact, in view of our computations in
Section \ref{s:convergence} below, it is important to introduce a fine decomposition of those terms.

\begin{lemma} \label{l:momentum}
Let $T>0$ be arbitrarily large, but fixed. For any $M\in\N$, define the quantities $\phi_{\veps,M}\,=\,S_M(\phi_\veps)$ and $\vc V_{\eps,M}\,=\,S_M(\vc V_\veps)$ as above.
\begin{enumerate}[(i)]
\item Both sequences $\big(\phi_{\eps,M}\big)_\eps$ and $\big(\vc V_{\eps,M}\big)_\eps$ are sequences of smooth functions in space, which are uniformly bounded (with respect to $\veps$, but \emph{not}
with respect to $M$) in the space $L^2_T(H^s)$, for any $s\geq0$ arbitrarily large.
 \item For any $M\in\N$ and any $\veps\in\,]0,1]$, we can write
\begin{equation} \label{eq:phi-decomp}
\phi_{\veps,M}\,=\,\vphi_{\veps,M}\,+\,\veps^{(2-\g)/\g}\,\pi_{\veps,M}\,,
\end{equation}
where we have the uniform estimates
\[
\sup_{M\in\N}\,\sup_{\veps\in\,]0,1]}\left\|\vphi_{\veps,M}\right\|_{L^\infty_T(L^2)\cap L^2_T(L^6)}\,\leq\,C\qquad \mbox{ and }\qquad
\forall\,s\geq0\,,\quad \sup_{\veps\in\,]0,1]} \left\|\pi_{\veps,M}\right\|_{L^\infty_T(H^s)}\,\leq\,C(s,M)\,,
\]
for suitable positive constants $C$ and $C(s,M)$, which depend only on $T>0$ and on the quantities in the brackets. In the case $\g\geq2$, one can simply take $\pi_{\veps,M}\equiv0$.
 \item For any $M\in\N$ and $\veps\in\,]0,1]$, one has 
\begin{equation} \label{eq:mom-decomp}
\vc V_{\eps,M}\,=\,\mc V_{\eps.M}\,+\,\eps\,\vc W_{\eps,M}\,,
\end{equation}
with the uniform bounds
\[
\sup_{M\in\N}\,\sup_{\veps\in\,]0,1]}\left\|\mc V_{\eps,M}\right\|_{L^\infty_T(L^2)}\,\leq\,C\qquad\quad \mbox{ and }\qquad\quad
\forall\,s\geq0\,,\quad \sup_{\veps\in\,]0,1]} \left\|\vc W_{\veps,M}\right\|_{L^2_T(H^s)}\,\leq\,C(s,M)\,.
\]
\item We can also write
\begin{equation} \label{eq:mom-decomp2}
\vc V_{\eps,M}\,=\,\wtilde{\mc V}_{\eps.M}\,+\,\eps\,\wtilde{\vc W}_{\eps,M}\,,
\end{equation}
with the uniform estimates
\[
\sup_{M\in\N}\,\sup_{\veps\in\,]0,1]}\left\|\wtilde{\mc V}_{\eps,M}\right\|_{L^2_T(W^{1,p_1})}\,\leq\,C\qquad\quad \mbox{ and }\qquad\quad
\forall\,s\geq0\,,\quad \sup_{\veps\in\,]0,1]} \left\|\wtilde{\vc W}_{\veps,M}\right\|_{L^2_T(H^s)}\,\leq\,C(s,M)\,.
\]
\end{enumerate}
\end{lemma}

\pf
The properties claimed in item (i) are an immediate consequence of the bounds for the families $\big(\phi_\veps\big)_\veps$ and $\big(\vc V_\veps\big)_\veps$, combined with
the smoothing effect of the low frequency cut-off operators $S_M$. Keep in mind \eqref{ub:phi_ess}, \eqref{ub:phi_res} and Propositions \ref{p:D-phi} and \ref{p:momentum}.

The decompositions of items (iii) and (iv), together with the corresponding uniform estimates, also follow from Proposition \ref{p:momentum}, simply setting 
$\mc V_{\eps,M}\,:=\,S_M(\mc V_\eps)$ and $\vc W_{\eps,M}\:=\,S_M(\vc W_\veps)$, and similarly for $\wtilde{\mc V}_{\veps,M}$ and $\wtilde{\vc W}_{\veps,M}$.

Finally, let us prove item (ii). We start by considering the case $1<\g<2$. In this case, we decompose
\[
\phi_\veps\,=\,\big[\phi_\veps\big]_\ess\,+\,\big[\phi_\veps\big]_\res\,.
\]
Thus, if we set $\vphi_{\veps,M}\,:=\,S_M\left(\big[\phi_\veps\big]_\ess\right)$, estimates \eqref{ub:phi_ess} and \eqref{est:rho-B} easily imply the uniform boundedness (\emph{both} with respect to $M\in\N$
and $\veps\in\,]0,1]$) of the sequence $\big(\vphi_{\veps,M}\big)_{\veps,M}$.
Next, we define
\[
\pi_\veps\,:=\,\frac{1}{\veps^{(2-\g)/\g}}\,\big[\phi_\veps\big]_\res\qquad\qquad\mbox{ and }\qquad\qquad \pi_{\veps,M}\,:=\,S_M(\pi_\veps)\,,
\]
and we conclude with the help of \eqref{ub:phi_res}. When $\g\geq2$, instead, one has that $\big(\phi_\veps\big)_\veps$ is uniformly bounded in $L^\infty_T(L^2)\cap L^2_T(H^1)$,
thanks to Proposition \ref{p:D-phi}. Hence, in this case one can simply define $\vphi_{\veps,M}\,=\,S_M(\phi_\veps)$ and $\pi_{\veps,M}\,=\,0$.
This completes the proof of the lemma. \qed

\section{Passage to the limit} \label{s:convergence}

In this section, we finish the proof of our main theorem, showing the convergence (up to an extraction) of weak solutions to \eqref{main} to weak solutions to the target system \eqref{eq:limit}.

The main problem is to pass to the limit in the most non-linear term, namely the convective term in the momentum equation. This convergence will be proved in the following two subsections,
by use of a \emph{compensated compactness} argument. 
Finally, in Subsection \ref{ss:final-comp} we will take care of the other terms, and complete the proof of the convergence.

Before going into the details, let us recall that convergence will be proved for any test function lying in the kernel of the singular perturbation operator, namely (in view of Proposition
\ref{p:constraints}) for any test function
\begin{equation} \label{eq:test}
\psi\,\in\,\mc D\big(\R_+\times\Omega;\R^3\big)\,,\qquad\qquad \mbox{ such that }\qquad 
\div\big(b\,\psi\big)\,=\,0\,.
\end{equation}

Also, it is useful to introduce the following notation: we denote by $R_{\eps, M}$ any remainder term, that is any term such that
\begin{equation} \label{eq:def-rem}
\lim_{M\ra+\infty}\,\limsup_{\eps\ra0^+} \left|  \int_0^T \langle R_{\eps, M}\,,\,\psi\rangle_{\mc D'\times \mc D}\,\dt \right|\,=\, 0\,,
\end{equation}
for some given time $T>0$ and test function $\psi\in\mc D\big([0,T[\,\times\Omega;\R^3\big)$ taken as above.
Similarly, we will use the notation $\mc R_{\ep,M}$ for any scalar term such that
\begin{equation} \label{eq:def-rem_int}
\lim_{M\ra+\infty}\,\limsup_{\eps\ra0^+}\; \mc R_{\ep,M}\,=\,0\,.
\end{equation}
Typically, we will have
\[
\mc R_{\veps,M}\,=\,\int_0^T \langle R_{\eps, M}\,,\,\psi\rangle_{\mc D'\times \mc D}\,\dt\,.
\]
In the next computations, the precise values of $R_{\veps,M}$ and $\mc R_{\veps,M}$ may change from one line to another.

\subsection{Approximation of the convective term}\label{ss:approx}

Passing to the limit in the convective term is based on a compensated compactness argument, following work \cite{Masm} of Masmoudi.
This argument relies on using \emph{algebraically} the structure of the wave system \eqref{eq:wave} and performing direct computations on it. Of course, for that, we need first of all
to smooth out all the quantities with respect to the space variables: this is the scope of the next lemma.

Notice that the approximation argument here is delicate, due to the degeneracy of the system close to vacuum. The cold pressure term $p_c$ will be of great help.
\begin{lemma} \label{l:convective}
For any $T>0$ fixed and any smooth test function $\psi$ as in \eqref{eq:test}, such that $\Supp\psi\subset[0,T[\,\times\Omega$, we have
\[
\lim_{M\ra+\infty}\,\limsup_{\eps\ra0^+}\left|\intTO{\vr_\eps\vu_\eps\otimes\vu_\eps:\nabla\psi}\,-\,\intTO{\frac{1}{b}\,\vc V_{\eps,M}\otimes \vc V_{\eps,M}:\nabla\psi}\,\right|\,=\,0\,.
\]
\end{lemma}

\pf
We start by observing that, in view of Propositions \ref{p:u} and \ref{p:momentum}, we have $\big(\vc V_\eps\big)_\eps\,\subset\,L^2_T(L^{3/2})$ and
$\big(\vu_\eps\big)_\eps\,\subset\,L^2_T(L^{p})$ for all $p\in[1,p_2]$, where $p_2\,=\,6\k/(\k+3)$. We also notice that, under our assumption $\k>3$, we have $2/3\,+\,1/p_2<1$, \tsl{i.e.} $p_2>3$.
Thus, given any $M\in\N$, we can write
\begin{align*}
I\,&:=\,\intTO{\vr_\eps\vu_\eps\otimes\vu_\eps:\nabla\psi}\, \\
&=\,\intTO{\frac1b\,\vc V_\eps\otimes S_M(b\,\vu_\eps):\nabla\psi}\,+\,
\intTO{\frac1b\,\vc V_\eps\otimes (\Id-S_M)(b\,\vu_\eps):\nabla\psi}\,, % \\
\end{align*}
where all terms are well-defined, and the last term on the right can be bounded as follows:
\[
\left|\intTO{\frac1b\,\vc V_\eps\otimes (\Id-S_M)(b\,\vu_\eps):\nabla\psi}\,\right|\,\lesssim\,\left\|\vc V_\eps\right\|_{L^2_T(L^{3/2})}\,\left\|(\Id-S_M)(b\,\vu_\eps)\right\|_{L^2_T(L^3)}\,.
\]
At this point, taking advantage of Bernstein's inequalities of Lemma \ref{l:bern} in the Appendix, we can bound
\begin{align}
 \left\|(\Id-S_M)(b\,\vu_\eps) \right\|_{L^{3}}\,&\leq\,\sum_{j\geq M}\left\|\Delta_j(b\,\vu_\eps)\right\|_{L^{3}}\,\lesssim\,\sum_{j\geq M}2^{-j}\left\|\Delta_j\nabla(b\,\vu_\eps)\right\|_{L^{3}}
\label{est:Id-S} \\
 &\lesssim\,\sum_{j\geq M}2^{j\left[-1\,+\,3\left(\frac{1}{p_1}-\frac{1}{3}\right)\right]}\left\|\Delta_j(b\,\vu_\eps)\right\|_{L^{p_1}}\,\leq\,
2^{-\alpha M}\,\left\|\nabla(b\,\vu_\eps)\right\|_{B^0_{p_1,\infty}}\,\sum_{j\geq 0}2^{-\alpha j}\,, \nonumber % \\
% &\leq\,C\,2^{-M/3}\,\left\|\nabla \vc V_\eps\right\|_{L^{9/8}}\,.
\end{align}
where $p_1=2\k/(\k+1)$ has been defined in Proposition \ref{p:u} and we have set $\alpha\,=\,(\k-3)/2\k\,>\,0$.
Owing to the embedding $L^p\,\hookrightarrow\,B^0_{p,\infty}$ for any $p\in[1,+\infty]$, we finally get
\[
 \left\|(\Id-S_M)(b\,\vu_\eps) \right\|_{L^{3}}\,\lesssim\,2^{-\alpha M}\,\left\|\nabla(b\,\vu_\eps)\right\|_{L^{p_1}}\,\lesssim\,2^{-\alpha M}\,.
\]
As a result, the previous computations show that, in the sense of \eqref{eq:def-rem_int}, one has
\[
I\,=\,%\int^T_0\intO{\vr_\eps\vu_\eps\otimes\vu_\eps:\nabla\psi}\,=\,
\int^T_0\intO{\frac1b\,\vc V_\eps\otimes S_M(b\,\vu_\eps):\nabla\psi}\,\dt\,+\,\mc R_{\veps,M}\,.
\]

Next, we use again the uniform bound $\big(\vu_\eps\big)_\eps\,\subset\,L^2_T(L^{p_2})$, together with the fact that $\big(S_M(b\,\vu_\eps)\big)_\eps$ is uniformly bounded (with respect to $\eps$, but
\emph{not} to $M$) in the space $L^2_T(L^\infty)$, to get
\begin{align} \label{eq:I_interm}
\hspace{-0.5cm}&\int^T_0\!\!\!\intO{\frac1b\,\vc V_\eps\otimes S_M(b\,\vu_\eps):\nabla\psi}\,\dt \\
&\qquad\qquad\qquad=\,\int^T_0\!\!\!\intO{\vu_\eps\otimes S_M(b\,\vu_\eps):\nabla\psi}\,\dt\,+\,
\veps\int^T_0\!\!\!\intO{\frac{1}{b}\,\phi_\eps\,\vu_\eps\otimes S_M(b\,\vu_\eps):\nabla\psi}\,\dt \nonumber \\
 &\qquad\qquad\qquad=\,\intTO{\frac{1}{b}\,S_M\big(b\,\vu_\eps\big)\otimes S_M(b\,\vu_\eps):\nabla\psi}\,+\,\mc R_{\veps,M}\,. \nonumber
\end{align}
Indeed, first note that, due to \eqref{ub:phi_higher}, we can write, for $\g<2$, the following estimate:
\[
\left|\intTO{\frac1b\,\phi_\eps\,\vu_\eps\otimes S_M(b\,\vu_\eps):\nabla\psi}\,\dt\right|\,\lesssim\,\left\|\phi_\eps\right\|_{L^\infty_T(L^{3\g/(3-\g)})}\,\left\|\vu_\veps\right\|_{L^2_T(L^{p_2})}\,
\left\|S_M(b\,\vu_\eps)\right\|_{L^2_T(L^\infty)}\,\lesssim\,1\,.
\]
We observe that this makes sense whenever $(3-\g)/3\g\,+\,(\k+3)/6\k\,\leq\,1$, hence for $\k\,\geq\,3\g/(7\g-6)$. But, for $\g>1$, one always has $3\g/(7\g-6)<3$, therefore the previous
estimate is satisfied for all $\kappa>3$. Notice that, since $p_2>3$, we can repeat the exact same computations also when $\g\geq2$, up to use the right bound from \eqref{ub:phi_higher}.

On the other hand, after noticing that $1/p_1\,+\,1/p_2\,<\,1$ for $\k> 3$ and arguing in a similar way as above, we can estimate
\begin{align*}
\left|\int^T_0\!\!\!\intO{\frac1b\,(\Id-S_M)(b\,\vu_\eps)\otimes S_M(b\,\vu_\eps):\nabla\psi}\,\dt\right|\,&\lesssim\,\left\|(\Id-S_M)(b\,\vu_\eps)\right\|_{L^2_T(L^{p_1})}\,
\left\|S_M(b\,\vu_\eps)\right\|_{L^2_T(L^{p_2})} \\
&\lesssim\,2^{- M}\,\left\|\nabla(b\,\vu_\eps)\right\|^2_{L^{p_1}}\,\lesssim\,2^{- M}\,.
\end{align*}
Thus, we have proven that the last equality in \eqref{eq:I_interm} holds true.

At this point, to conclude the argument we may use the decomposition of Lemma \ref{l:momentum}. Indeed, keeping in mind the definition of the vector fields $\wtilde{\mc V}_\ep$
given in Proposition \ref{p:momentum} above, we notice that
$S_M(b\,\vu_\ep)\,=\,\wtilde{\mc V}_{\ep,M}\,=\,\vc V_{\ep,M}\,-\,\ep\,\wtilde{\vc W}_{\ep,M}$. Thus, using the bounds collected in item (iv) of Lemma \ref{l:momentum}
we easily see that
\[
\intTO{\frac{1}{b}\,S_M\big(b\,\vu_\eps\big)\otimes S_M(b\,\vu_\eps):\nabla\psi}\,=\,\intTO{\frac{1}{b}\,\vc V_{\veps,M}\otimes \vc V_{\veps,M}:\nabla\psi}\,+\,\mc R_{\veps,M}\,.
\]
This last equality finally ends the proof of the lemma. \qed

\subsection{Compensated compactness} \label{ss:comp-comp}

Owing to Lemma \ref{l:convective}, it is enough to compute the limit of the approximate convective term
\[
-\,\intTO{\frac{1}{b}\,\vc V_{\eps,M}\otimes \vc V_{\eps,M}:\nabla\psi}\,=\,\intTO{\div\left(\frac{1}{b}\,\vc V_{\eps,M}\otimes \vc V_{\eps,M}\right)\cdot\psi}\,,
\]
for any test function $\psi\in\mc D\big([0,T[\,\times\Omega;\R^3\big)$ belonging to the kernel of the singular perturbation operator, namely such that $\div(b\psi)\,=\,0$.
Here above, we have performed an integration by parts, because the vector fields $\vc V_{\eps,M}$ are smooth in the space variable.

\subsubsection{Preliminary reductions} \label{sss:prelim}
Now, we compute
\begin{align}
\div\left(\frac{1}{b}\,\vc V_{\eps,M}\otimes \vc V_{\eps,M}\right)\,&=\,\frac{1}{b}\,\Div(\vc V_{\ep,M})\, \vc V_{\ep,M}\,+\,\vc V_{\eps,M}\cdot\nabla\left(\frac{1}{b}\,\vc V_{\veps,M}\right) \label{eq:comp-conv} \\
&=\,\frac{1}{b}\,\Div (\vc V_{\ep,M})\, \vc V_{\ep,M}\,+\,\frac{b}{2}\,\Grad\left|\frac{1}{b}\,\vc V_{\veps,M}\right|^2\,-\,\curl\left(\frac{1}{b}\,\vc V_{\eps,M}\right)\times \vc V_{\eps,M}\,, \nonumber
\end{align}
where the symbol $\times$ denotes the ususal external product of vectors in $\R^3$ and, for a $3$-D vector field $\vc U$, we have $\curl \vc U:= \nabla\times \vc U$.
Notice that the second term in the last line identically vanishes, whenever tested against a test function $\psi$ as in \eqref{eq:test}. %$\psi$ such that $\div\big(b\,\psi\big)\,=\,0$.
Thus, this term contributes as a remainder $R_{\veps,M}$ to the limit, in the sense of relation \eqref{eq:def-rem}.

Therefore, resorting to the first equation in \eqref{eq:wave_M} for dealing with the first term, we can write
\begin{align*}
 \div\left(\frac{1}{b}\,\vc V_{\eps,M}\otimes \vc V_{\eps,M}\right)\,&=\,-\,\frac{\veps}{b}\,\d_t\phi_{\eps,M}\,\vc V_{\eps,M}\,-\,\curl\left(\frac{1}{b}\,\vc V_{\eps,M}\right)\times \vc V_{\eps,M}\,+\,R_{\eps,M} \\
 &=\,\frac{\veps}{b}\,\phi_{\eps,M}\,\d_t\vc V_{\eps,M}\,-\,\curl\left(\frac{1}{b}\,\vc V_{\eps,M}\right)\times \vc V_{\eps,M}\,+\,R_{\eps,M}\,,
\end{align*}
where, in the second step, we have included the total time derivative $\veps\,\d_t\big(\phi_{\eps,M}\,\vc V_{\eps,M}\big)$ into the remainder $R_{\eps,M}$. Indeed, the time derivative
can be put on the test function and the family $\big(\phi_{\eps,M}\,\vc V_{\eps,M}\big)_{\eps}$ is uniformly bounded in \tsl{e.g.} $L^2_T(L^{2})$, owing to item (i) of Lemma \ref{l:momentum}.

At this point, we use the second equation in \eqref{eq:wave_M} to deal with the term presenting the time derivative, and we get
\begin{align*}
&\div\left(\frac{1}{b}\,\vc V_{\eps,M}\otimes \vc V_{\eps,M}\right) \\
&\qquad\qquad=\,-\,\phi_{\eps,M}\,\Grad\big(H''(b)\,\phi_{\veps,M}\big)\,+\,\frac{\veps}{b}\,\phi_{\eps,M}\,\vc f_{\eps,M}\,+\,
\frac{1}{b}\,\phi_{\eps,M}\,\vc h_{\eps,M}\,-\,\curl\left(\frac{1}{b}\,\vc V_{\eps,M}\right)\times \vc V_{\eps,M}\,+\,R_{\eps,M}\,.
\end{align*}
Owing to \eqref{est:f} and item (i) of Lemma \ref{l:momentum}, it is clear that the second term in the right-hand side is a remainder, in the sense of \eqref{eq:def-rem}.
Next, we claim that also the third term, \tsl{i.e.} $\big(\phi_{\eps,M}\,\vc h_{\eps,M}/b\big)_{\veps,M}$, is a remainder. For proving this claim, in the case $\g\geq2$ it is enough
to employ Proposition \ref{p:D-phi} and recall that estimate \eqref{est:h} holds true with the Lebesgue exponent $\g$ replaced by $2$.
When $\g\in\,]1,2[\,$, instead, we notice that the piece of information coming
from \eqref{ub:phi_higher} is not good enough to cover all possible values of $\g$ in that interval. Our approach is instead based on the use of item (ii) of Lemma \ref{l:momentum},
which allows us to write
\[
\frac{1}{b}\,\phi_{\eps,M}\,\vc h_{\eps,M}\,=\,\frac{1}{b}\,\vphi_{\eps,M}\,\vc h_{\eps,M}\,+\,\veps^{(2-\g)/\g}\,\frac{1}{b}\,\pi_{\eps,M}\,\vc h_{\eps,M}\,.
\]
Thanks to \eqref{est:h} and item (ii) of Lemma \ref{l:momentum},  we can estimate
\begin{align*}
\left\|\frac{1}{b}\,\vphi_{\eps,M}\,\vc h_{\eps,M}\right\|_{L^1_T(L^{6\g/(6-\g)})}\,\lesssim\,\left\|\vphi_{\veps,M}\right\|_{L^2_T(L^6)}\,\left\|\vc h_{\veps,M}\right\|_{L^2_T(L^\g)}\,\lesssim\,2^{-M}\,,
\end{align*}
which implies that this term satisfies \eqref{eq:def-rem}. In addition, it is easy to see that also $\big(\veps^{(2-\g)/\g}\,\pi_{\eps,M}\,\vc h_{\eps,M}/b\big)_\veps$ verifies \eqref{eq:def-rem}.
Indeed, recall that  in the case $\g\geq2$, one can simply take $\pi_{\veps,M}\equiv0$. For $1<\g<2$ one has only to notice that
\begin{align*}
\left\|\veps^{(2-\g)/\g}\,\frac{1}{b}\,\pi_{\eps,M}\,\vc  h_{\eps,M}\right\|_{L^2_T(L^\g)}\,\lesssim\,\veps^{(2-\g)/\g}\,\left\|\pi_{\veps,M}\right\|_{L^\infty_T(L^\infty)}\,
\left\|\vc h_{\veps,M}\right\|_{L^\infty_T(L^\g)}\,\lesssim\,C(M)\,\veps^{(2-\g)/\g}\,.
\end{align*}

As a result of the previous computations, we infer that
\begin{equation} \label{eq:convect-1}
\div\left(\frac{1}{b}\,\vc V_{\eps,M}\otimes \vc V_{\eps,M}\right)\,=\,-\,\phi_{\eps,M}\,\Grad\big(H''(b)\,\phi_{\veps,M}\big)\,-\,\curl\left(\frac{1}{b}\,\vc V_{\eps,M}\right)\times \vc V_{\eps,M}\,+\,R_{\eps,M}\,.
\end{equation}

\subsubsection{Compactness of the rotational part} \label{sss:curl-cpt}
The next lemma takes care of the convergence of the $\curl$ term in \eqref{eq:convect-1}.

\begin{lemma} \label{l:curl}
Denote by $\vc V$ the weak-limit of $\big(\vc V_\veps\big)_\veps$ identified in Corollary \ref{c:momentum}, so that $\vc V_M\,=\,S_M\vc V$ is the weak-limit
of $\big(\vc V_{\veps,M}\big)_\veps$ when $\veps\ra0^+$.
Then, for any $\psi\in\mc D\big(\R_+\times\Omega;\R^3\big)$, one has
\[
\lim_{M\ra+\infty}\limsup_{\veps\ra0^+}\left(\int^T_0\!\!\!\intO{\curl\left(\frac{1}{b}\,\vc V_{\eps,M}\right)\times \vc V_{\eps,M}\cdot\psi},\dt\,-\,
\int^T_0\!\!\!\intO{S_M\curl\left(\frac{1}{b}\,\vc V\right)\times \vc V_M\cdot\psi}\,\dt\right)\,=\,0\,.
\]
\end{lemma}
\pf We start the proof by recalling that, owing to Proposition \ref{p:momentum}, we have $\big(\vc V_\eps\big)_\eps\,\subset\,L^2_T(W^{1,p_3})$,
hence $\big(\curl \vc V_\eps\big)_\eps\,\subset\,L^2_T(L^{p_3})$. Observe that, by dual Sobolev embeddings (see \tsl{e.g.} Theorem 0.5 of \cite{F-N}),
we have that, for any compact $K\subset \Omega$, the space $L^{p_3}(K)$ is compactly embedded into $H^{-2}(K)$, for instance.

Next, consider the (not regularised) wave system \eqref{eq:wave}. Dividing the momentum equation by $b$ and then taking the $\curl$
(which simply consists in an adequate choice of the test function in the weak formulation), we deduce 
\[
\pt\curl\lr{\frac{1}{b} \vc V_\ep} =\curl\lr{\frac{\vc f_\ep}{b}}\,.
\]
In turn, this relation, together with \eqref{unif-b:f} and the fact that $b\in C^3(\Omega)$,
implies that $\big(\d_t\curl\lr{\frac{1}{b} \vc V_\ep}\big)_\ep$ is uniformly bounded in the space $L^2_T(H^{-s})$, for any $s>7/2$.

Putting together these pieces of information and applying the Aubin-Lions lemma (see \tsl{e.g.} Lemma 3.7 of \cite{Tsai}), we gather that
\[
\bigg(\curl\lr{\frac{1}{b} \vc V_\ep}\bigg)_\ep\qquad\qquad\mbox{ is compact in }\qquad\qquad L^2_T\big(H^{-2}(K)\big)\,,
\]
for any compact set $K\subset\Omega$.
This implies that, for any $M\in\N$ and any compact $K\subset\Omega$ fixed, the sequence
\begin{equation} \label{eq:curl-cpt}
\bigg(S_M\curl\lr{\frac{1}{b} \vc V_\ep}\bigg)_\ep\qquad\qquad\mbox{ is compact (with respect to $\veps$) in }\qquad\qquad L^2_T\big(L^{2}(K)\big)\,.
\end{equation}

Next, we write
\begin{align*}
&\intTO{\curl\left(\frac{1}{b}\,\vc V_{\eps,M}\right)\times \vc V_{\eps,M}\cdot\psi} \\
&\qquad =\,\intTO{S_M\curl\left(\frac{1}{b}\,\vc V_{\eps}\right)\times \vc V_{\eps,M}\cdot\psi}\,+\,
\intTO{\curl\left(\left[\frac{1}{b},S_M\right]\,\vc V_{\eps}\right)\times \vc V_{\eps,M}\cdot\psi}\,.
\end{align*}
By what we have just said, we have that
\[
\lim_{\veps\ra0^+}\intTO{S_M\curl\left(\frac{1}{b}\,\vc V_{\eps}\right)\times \vc V_{\eps,M}\cdot\psi}\,=\,
\intTO{S_M\curl\left(\frac{1}{b}\,\vc V\right)\times \vc V_{M}\cdot\psi}\,.
\]
On the other hand, using the embedding $\big(\vc V_\eps\big)_\eps\,\subset\,L^2_T(W^{1,p_3})$ again, we have
\[
 \left\|\curl\left(\left[\frac{1}{b},S_M\right]\,\vc V_{\eps}\right)\right\|_{L^2_T(L^{p_3})}\,\lesssim\,2^{-M}\,,
\]
whereas item (iv) of Lemma \ref{l:momentum} and Bernstein's inequalities (see Lemma \ref{l:bern} in the Appendix) imply
\[
\left\|\vc V_{\veps,M}\right\|_{L^2_T(L^{p_3'})}\,\lesssim\,2^{3M(\frac{1}{p_2}\,-\,\frac{1}{p_3'})}\,+\,\veps\,C_M\,,
\]
for a constant $C_M$ which blows up when $M\ra+\infty$, but which is uniform in $\veps>0$. Observe that
\[
3\left(\frac{1}{p_2}\,-\,\frac{1}{p_3'}\right)\,=\,3\left(\frac{\k+3}{6\k}\,-\,1\,+\,\frac{5\k+3}{6\k}\right)\,=\frac{3}{\kappa}\,,\qquad\qquad \mbox{ with }\qquad \frac{3}{\k}\,<\,1\,.
\]
Thus, we deduce that
\[
\lim_{M\ra+\infty}\limsup_{\veps\ra0^+}\intTO{\curl\left(\left[\frac{1}{b},S_M\right]\,\vc V_{\eps}\right)\times \vc V_{\eps,M}\cdot\psi}\,=\,0\,.
\]

The proof of the lemma is now completed. \qed

\subsubsection{Handling the pressure term} \label{sss:pressure}
Before computing the limit with respect to $M\ra+\infty$ in the previous lemma, let us treat the first term in the right-hand side of \eqref{eq:convect-1}.
As a matter of fact, we need to couple it with a term coming from the pressure function, namely 
\[ %\begin{equation} \label{eq:rem-pressure}
-\,\frac{1}{\veps^2}\,\intO{\Pi(\vr_\veps;b)\,\Div\psi}\,,
\] %\end{equation}
where $\Pi(\vr_\eps;b)$ has been defined in Lemma \ref{l:pressure}. We will use a fundamental cancellation (appearing after regularisation), which is already present in \cite{Masm}.
For this, we need the following preparatory lemma.

\begin{lemma} \label{l:Pi}
For any $T>0$ fixed and any test function $\psi$ as in \eqref{eq:test}, such that $\Supp\psi\subset[0,T[\,\times\Omega$, the following relation holds true, in the sense of \eqref{eq:def-rem_int}:
\[
-\frac{1}{\veps^2}\,\intTO{\Pi(\vr_\veps;b)\,\Div\psi}\,=\,-\,\frac{1}{2}\int^T_0\!\!\!\int_{\Omega^\ep_\ess}{p''(b)\,\phi_\ep^2\,\Div\psi}\,\dx\,\dt\,+\,\mc R_{\veps,M}\,.
\]
\end{lemma}

\pf
Let us start by decomposing the term on the left-hand side of the claimed equality into two integrals, one on the essential set and the other on the residual set:
\[
\frac{1}{\veps^2}\,\intTO{\Pi(\vr_\veps;b)\,\Div\psi}\,=\,\frac{1}{\veps^2}\,\int^T_0\!\!\!\int_{\Omega^\veps_\ess}{\Pi(\vr_\veps;b)\,\Div\psi}\,\dx\,\dt\,+\,
\frac{1}{\veps^2}\,\int^T_0\!\!\!\int_{\Omega^\veps_\res}{\Pi(\vr_\veps;b)\,\Div\psi}\,\dx\,\dt\,.
\]
Notice that, in order to treat the integral over $\Omega^\veps_\res$, the bounds of Lemma \ref{l:pressure} are not enough.
Instead, in order to take advantage of the smallness bounds in \eqref{ub:dens-res}, we need to introduce a finer decomposition and split that term further into two pieces.

To begin with, let us consider $\Pi_E(\vre;b)\,=\,p_E(\vr_\eps)-p_E(b)-p_E'(b)(\vre-b)$ only.
On the one hand, thanks the uniform bounds of \eqref{est:rho-B}, we have, for some $z_\ep=z_\ep(t,x)$ belonging to the interval $\,]\vr_\ep(t,x),b(x)[\,$, the estimate
\begin{align*}
\left|\frac{1}{\veps^2}\,\int_{\Omega^\eps_{\res,B}}{\Pi_E(\vr_\veps;b)\,\Div\psi}\,\dx\right|\,&\lesssim\,
\frac{1}{\veps}\,\int_{\Omega^\eps_{\res,B}}{\Big(p'_E(z_\ep)+p'_E(b)\Big)\,\big|\phi_\ep\big|\,\big|\Div\psi\big|}\,\dx \\
&\lesssim\,\frac{1}{\veps}\,\left\|\mathds{1}_{\Omega^\veps_{\res,B}}\,\phi_\ep\right\|_{L^2_T(L^6)}\,\Big(\mc L(\Omega^\veps_\res)\Big)^{5/6}\,\lesssim\,\veps^{-1+5/3}\,=\,\veps^{2/3}\,,
\end{align*}
which tells us that the contribution coming from this term is a remainder, in the sense of \eqref{eq:def-rem}. On the other hand, % For the other term, we proceed similarly:
proceeding similarly and employing the assumption $p_E'(z)\approx z^{\g-1}$, we get
\begin{align*}
\left|\frac{1}{\veps^2}\,\int_{\Omega^\ep_{\res,UB}}{\Pi_E(\vr_\veps;b)\,\Div\psi}\,\dx\right|\,&\lesssim\,
\frac{1}{\veps}\,\int_{\Omega^\ep_{\res,UB}}{\vr_\ep^{\g-1}\,\big|\phi_\ep\big|\,\big|\Div\psi\big|}\,\dx\,.% \\
\end{align*}
Now, owing to \eqref{est:rho-UB}, we gather the uniform bound
\begin{equation} \label{est:rho-ep}
\left(\frac{1}{\veps}\,\mathds{1}_{\Omega^\veps_{\res,UB}}\,\vr_\veps^{\g-1}\right)_\veps\,\subset\,L^{r}_T(L^{p})\,,\qquad\qquad 
r\,:=\,\frac{\g}{\g-1}\quad\mbox{ and }\quad p\,:=\,\frac{3\,\g}{\g-1}\,,
\end{equation}
whereas we can use \eqref{ub:phi_higher} to control the $\phi_\veps$ term. After checking that
\[
\forall\;1\,<\,\g\,\leq\,2\,,\qquad\qquad\frac{\g-1}{\g}\,+\,\frac{1}{2}\,\leq\,1\qquad\mbox{ and }\qquad \frac{1}{q}\,:=\,\frac{3-\g}{3\,\g}\,+\,\frac{\g-1}{3\,\g}\,<\,1\,,
\]
we can thus bound
\begin{align*}
\left|\frac{1}{\veps^2}\,\int^T_0\!\!\!\int_{\Omega^\ep_{\res,UB}}{\Pi_E(\vr_\veps;b)\,\Div\psi}\,\dx\,\dt\right|\,&\lesssim\,
\left(\frac{1}{\veps}\,\left\|\mathds{1}_{\Omega^\veps_{\res,UB}}\,\vr_\eps^{\g-1}\right\|_{L^{r}_T(L^p)}\right)\,\left\|\phi_\ep\right\|_{L^2_T(L^{3\g/(3-\g)})}\,\Big(\mc L(\Omega_\res)\Big)^{1-1/q} \\
&\lesssim\,\veps^{2(1-1/q)}\,.
\end{align*}
A direct computation shows that $1-1/q\,>\,0$. Therefore, for $1<\g\leq2$, we have proved that the contribution coming from the integral over the residual set is indeed a remainder, in the sense of \eqref{eq:def-rem}.

In the case $\g>2$, instead, we interpolate between the $L^\infty_T(L^{\g/(\g-1)})$ bound on $\big(\mathds{1}_{\Omega^\veps_{\res,UB}}\,\vr_\eps^{\g-1}\big)_\veps$, which comes from
\eqref{ub:dens-res}, and the bounds provided by \eqref{est:rho-ep}. We find that
\[
\left(\mathds{1}_{\Omega^\veps_{\res,UB}}\,\vr_\eps^{\g-1}\right)_\veps\,\subset\,L^2_T(L^m)\,,\qquad\qquad m\,:=\,\frac{3\g}{2\g-3}\,\geq\,\frac65\,.
\]
Hence, owing to \eqref{ub:phi_higher} again, when $\g>2$ we can set $\theta\,:=\,\frac{\g-2}{2(\g-1)}$ and estimate
\begin{align*}
\left|\frac{1}{\veps^2}\,\int^T_0\!\!\!\int_{\Omega^\ep_{\res,UB}}{\Pi_E(\vr_\veps;b)\,\Div\psi}\,\dx\,\dt\right|\,&\lesssim\,
\frac{1}{\veps}\,\left\|\mathds{1}_{\Omega^\veps_{\res,UB}}\,\vr_\eps^{\g-1}\right\|_{L^{2}_T(L^m)}\,\left\|\phi_\ep\right\|_{L^2_T(L^{6})}\,\Big(\mc L(\Omega^\veps_\res)\Big)^{5/6-1/m} \\
&\lesssim\,\frac{1}{\veps}\,\left\|\mathds{1}_{\Omega^\veps_{\res,UB}}\,\vr_\eps^{\g-1}\right\|^\theta_{L^{\infty}_T(L^{\g/(\g-1)})}\,
\left\|\mathds{1}_{\Omega^\veps_{\res,UB}}\,\vr_\eps^{\g-1}\right\|^{1-\theta}_{L^{r}_T(L^{p})}\,\veps^{2(5/6-1/m)} \\
&\lesssim\,\veps^{2/3}\;\veps^{-2/m}\;\veps^{2\theta(\g-1)/\g}\;\veps^{1-\theta}\,=\,\veps^{4/3}\;\veps^{-(\g-2)/2(\g-1)}\,.
\end{align*}
A simple computation shows that
\[
\forall\;\g\,>\,\frac{2}{5}\,,\qquad\qquad\frac{4}{3}\,-\,\frac{\g-2}{2(\g-1)}\,>\,0\,,
\]
which finally implies that, also in the case $\g>2$, the integral over the residual set is a remainder, in the sense of \eqref{eq:def-rem}.

In the end, we have proved that, for any value of $\g>1$, the contribution coming from the integral of $\Pi_E$ over the residual set is a remainder, in the sense of relation \eqref{eq:def-rem}.

\medbreak
Next, let us consider the integral involving $\Pi_c(\vre;b)\,:=\,p_c(\vre)-p_c(b)-p'_c(b)(\vre-b)$ over the residual set. For this term, the roles of $\Omega^\veps_{\res,B}$ and
$\Omega^\veps_{\res,UB}$ are inverted. For instance, by Taylor formula at the first order we can write, for suitable $z_\veps(t,x)\in\,]b,\vre(t,x)[\,$, the following estimate:
\begin{align*}
\left|\frac{1}{\veps^2}\,\int_{\Omega^\eps_{\res,UB}}{\Pi_c(\vr_\veps;b)\,\Div\psi}\,\dx\right|\,&\lesssim\,
\frac{1}{\veps^2}\,\int_{\Omega^\eps_{\res,UB}}{\Big(p'_c(z_\ep)+p'_c(b)\Big)\,\vr_\ep\,\big|\Div\psi\big|}\,\dx \\
&\lesssim\,\frac{1}{\veps^2}\,\left\|\mathds{1}_{\Omega^\veps_{\res,UB}}\,\vr_\ep\right\|_{L^\g_T(L^{3\g})}\,\Big(\mc L(\Omega^\veps_\res)\Big)^{1-1/(3\g)}\,\lesssim\,\veps^{1-2/(3\g)}\,,
\end{align*}
which obviously converges to $0$ when $\veps\ra0^+$. In the previous computation, we have used \eqref{est:rho-UB} and \eqref{ub:dens-res} to absorb the negative powers of $\veps$.

In $\Omega_{\res,B}$, instead, Taylor formula and the fact that $z_\veps(t,x)\in\,]\vre(t,x),b[\,$ yield
\begin{align*}
&\hspace{-1cm}\left|\frac{1}{\veps^2}\,\int^T_0\!\!\!\int_{\Omega^\eps_{\res,B}}{\Pi_c(\vr_\veps;b)\,\Div\psi}\,\dx\,\dt\right| \\
&\qquad\lesssim\,
\frac{1}{\veps^2}\,\int^T_0\!\!\!\int_{\Omega^\eps_{\res,B}}{\Big(p'_c(z_\ep)+p'_c(b)\Big)\,(b-\vr_\ep)\,\big|\Div\psi\big|}\,\dx\,\dt\;\lesssim\;
\frac{1}{\veps^2}\,\int^T_0\!\!\!\int_{\Omega^\eps_{\res,B}}{\vr_\veps^{-\k-1}}\,\dx\,\dt\\
&\qquad\lesssim\,%\frac{1}{\veps^2}\,\int_{\Omega^\eps_{\res,B}}{\vr_\veps^{-\k-1}}\;\lesssim\;
\frac{1}{\veps^2}\,\left\|\mathds{1}_{\Omega^\eps_{\res,B}}\,\vre^{-\k/2}\right\|^2_{L^2_T(L^6)}\;\left\|\big[\vre^{-1}\big]_\res\right\|_{L^\infty_T(L^\k)}\;\Big(\mc L(\Omega^\veps_\res)\Big)^{2/3-1/\k}\,,
\end{align*}
which, in view of \eqref{est:rho-1} and \eqref{ub:dens-res}, is of order $O(\veps^{4/3})$. 

To sum up, we have shown that also the contribution from $\Pi_c$ over the residual set is a remainder, in the sense of \eqref{eq:def-rem}.

\medbreak
In light of what we have shown above, in order to complete the proof of the lemma 
it remains us to deal with the integral on the essential set. For this, we remark that, by Taylor formula, we can write
\begin{align*}
\frac{1}{\veps^2}\,\int^T_0\!\!\!\int_{\Omega^\ep_\ess}{\Pi(\vr_\veps;b)\,\Div\psi}\,\dx\,\dt\,=\,\frac{1}{2}\int^T_0\!\!\!\int_{\Omega^\ep_\ess}{p''(b)\,\phi_\ep^2\,\Div\psi}\,\dx\,\dt\,+\,
\frac{\veps}{6}\,\int^T_0\!\!\!\int_{\Omega^\ep_\ess}{p'''(z_\ep)\,\phi_\ep^3\,\Div\psi}\,\dx\,\dt\,,
\end{align*}
for some $z_\eps=z_\eps(t,x)$ belonging to the interval joining $b(x)$ and $\vr_\eps(t,x)$. At this point, we notice that the second integral on the right-hand side of the previous equality
can be bounded, with the help of \eqref{ub:phi_ess} and \eqref{est:rho-B}, in the following way:
\[
\left|\int^T_0\!\!\!\int_{\Omega_\ess}{p'''(z_\ep)\,\phi_\ep^{3}\,\Div\psi}\,\dx\,\dt\right|\,\lesssim\,\left\|\big[\phi\big]_\ess^2\right\|_{L^1_T(L^3)}\,
\left\|\big[\phi_\veps\big]_\ess\right\|_{L^\infty_T(L^2)}\,\lesssim\,1\,.
\]

Putting everything together, we finally get the claimed relation.
\qed

Of course, in the equality of Lemma \ref{l:Pi} there is no dependence of the remainder on the approximation parameter $M$.
However, in order for this information to be useful, we need to introduce the regularisation in the first term of the right-hand side. This is easy, thanks to Proposition \ref{p:D-phi} and
the fact that, on the essential set, the function $\phi_\veps$ possesses higher integrability in space (keep in mind property \eqref{est:rho-B}).

Indeed, for any $M\in\N$, we can write
\[
\big[\phi_\veps\big]_\ess^2\,=\,\big[\phi_\veps\big]_\ess\,\phi_{\veps,M}\,+\,\big[\phi_\veps\big]_\ess\,\big(\Id-S_M\big)\phi_\veps\,.
\]
Concerning the last term on the right, for any $T>0$ fixed, we proceed as follows: we use the uniform bounds $\big(\big[\phi_\veps\big]_\ess\big)_\veps\,\subset\,L^2_T(L^6)$ and, for $6/5\leq\g<2$,
the fact that
\begin{equation} \label{est:D-phi_e-M}
\left\|\big(\Id-S_M\big)\phi_\veps\right\|_{L^2_T(L^{\g})}\,\lesssim\,2^{-M}\,\left\|\big(\Id-S_M\big)\nabla\phi_\veps\right\|_{L^2_T(L^{\g})}\,\lesssim\,2^{-M}\,.
\end{equation}
This control follows from Bernstein's inequality (see Lemma \ref{l:bern} in the Appendix) and the uniform bounds of Proposition \ref{p:D-phi}, by arguing in a similar way as done in \eqref{est:Id-S}.
When $1<\g<6/5$, instead, we can use an interpolation argument between Lebesgue norms (because the inequality $3\g/(3-\g)>6/5$ holds for any $\g>6/7$) to get,
for a suitable $\theta\in\,]0,1[\,$, the following series of inequalities:
\begin{align*}
\left\|\big(\Id-S_M\big)\phi_\veps\right\|_{L^2_T(L^{6/5})}\,&\lesssim\,\left\|\big(\Id-S_M\big)\phi_\veps\right\|^\theta_{L^2_T(L^{\g})}\,
\left\|\big(\Id-S_M\big)\phi_\veps\right\|_{L^2_T(L^{3\g/(3-\g)})}^{1-\theta} \\
&\lesssim\,2^{-\theta M}\,\left\|\nabla\phi_\veps\right\|^\theta_{L^2_T(L^{\g})}\,\left\|\nabla\phi_\veps\right\|^{1-\theta}_{L^2_T(L^{\g})}\,\lesssim\,2^{-\theta M}\,.
\end{align*}
Finally, when $\g\geq2$, we can simply use Proposition \ref{p:D-phi} and relation \eqref{eq:LP-Sob} of the Appendix to find that $\left\|\big(\Id-S_M\big)\phi_\veps\right\|_{L^2_T(L^{2})}\,\lesssim\,2^{-M}$.

Next, we further decompose
\[
\big[\phi_\veps\big]_\ess\,\phi_{\veps,M}\,=\,\big[\phi_{\veps,M}\big]_\ess^2\,+\,\mathds{1}_{\Omega^\veps_\ess}\,\big(\Id-S_M\big)\phi_\veps\,\phi_{\veps,M}\,.
\]
At this point, we notice that all the uniform bounds satisfied by $\big(\phi_\ep\big)_\ep$ are also satisfied by $\big(\phi_{\ep,M}\big)_\ep$, uniformly with respect to
\emph{both} $M\in\N$ and $\ep\in\,]0,1]$. Thus, the same argument as above yields
\[
\left\|\mathds{1}_{\Omega^\veps_\ess}\,\big(\Id-S_M\big)\phi_\veps\,\phi_{\veps,M}\right\|_{L^1_T(L^{1})}\,\lesssim\,2^{-\theta M}\,,
\]
where $\theta=1$ if $\g\geq 6/5$, whereas $\theta\in\,]0,1[\,$ is as above in the case $1<\g<6/5$.

In the end, putting everything together, we have discovered that
\begin{equation} \label{eq:pressure-reg}
-\frac{1}{\veps^2}\,\int^T_0\!\!\!\intO{\Pi(\vr_\veps;b)\,\Div\psi}\,\dt\,=\,-\,\frac{1}{2}\int^T_0\!\!\!\int_{\Omega^\ep_\ess}{p''(b)\,\phi_{\ep,M}^2\,\Div\psi}\,\dx\,\dt\,+\,\mc R_{\veps,M}\,.
\end{equation}

\subsubsection{Coupling the pressure with the convective term}

We are now ready to deal with the first term appearing in the right-hand side of \eqref{eq:convect-1}. We have to pay attention, because here the signs are important.
We start by writing, for $\psi$ as in \eqref{eq:test},
\begin{align*}
-\intO{\phi_{\ep,M} \Grad\big(H''(b)\,\phi_{\ep,M}\big)\cdot\psi}\,&=\,-\,\frac{1}{2}\intO{H''(b)\,\Grad\phi_{\ep,M}^2\cdot\psi}
-\,\intO{\phi_{\ep,M}^2\, H'''(b)\,\psi\cdot\Grad b} \\
&=\,-\,\frac{1}{2}\intO{H''(b)\,\Grad\phi_{\ep,M}^2\cdot\psi}
+\,\intO{\phi_{\ep,M}^2\, H'''(b)\,b\,\,\Div\psi}\,, % \\
\end{align*}
where we have also used the fact that $\Div(b\,\psi)=0$. An integration by parts shows that
\begin{align*}
-\,\frac{1}{2}\intO{H''(b)\,\Grad\phi_{\ep,M}^2\cdot\psi}\,&=\,\frac{1}{2}\intO{\phi_{\veps,M}^2\,H''(b)\,\div\psi}\,+\,
\frac{1}{2}\intO{\phi_{\veps,M}^2\,H'''(b)\,\nabla b\cdot\psi} \\
&=\,\frac{1}{2}\intO{\phi_{\veps,M}^2\,H''(b)\,\div\psi}\,-\,\frac{1}{2}\intO{\phi_{\veps,M}^2\,H'''(b)\, b\,\div\psi}\,.
\end{align*}
At this point, we insert this expression into the previous one; after exploiting the definition $H''(z)\,=\,p'(z)/z$ for all $z>0$, we finally gather
\[
-\intO{\phi_{\ep,M} \Grad\big(H''(b)\,\phi_{\ep,M}\big)\cdot\psi}\,=\,\frac{1}{2}\,\intO{\phi_{\ep,M}^2\,\Div\psi\,p''(b)}\,.
\]
Owing to item (i) of Lemma \ref{l:momentum} and estimate \eqref{ub:residual-set}, we easily see that the integral over the residual set is small.
Therefore, after integrating also in time, in the end we get
\begin{equation*}
 -\intTO{\phi_{\ep,M} \Grad\big(H''(b)\,\phi_{\ep,M}\big)\cdot\psi}\,=\,\frac{1}{2}\,\int^T_0\!\!\!\int_{\Omega^\ep_\ess}{\phi_{\ep,M}^2\,\Div\psi\,p''(b)}\,\dx\,\dt\,+\,\mc R_{\eps,M}\,.
\end{equation*}
The fundamental point, here, is that the first term on the right-hand side exactly cancels out with the term coming from \eqref{eq:pressure-reg}.

\subsubsection{Limit of the convective term: conclusion} \label{sss:conv-concl}

Putting together Lemma \ref{l:convective} and the computations of Subsection \ref{ss:comp-comp}, we finally discover that,
for any test function $\psi$ belonging to the kernel of the singular perturbation operator, namely such that \eqref{eq:test} holds true, one has
\begin{align*}
&\lim_{\eps\ra0^+}\left(\intTO{\vr_\eps\vu_\eps\otimes\vu_\eps:\nabla\psi}\,-\,\frac{1}{\veps^2}\,\intTO{\Pi(\vr_\veps;b)\,\Div\psi}\,\dt\right) \\
&\qquad\qquad\qquad\qquad\qquad\qquad\qquad\qquad\qquad\qquad =\,\lim_{M\ra+\infty}\intTO{S_M\curl\left(\frac{1}{b}\,\vc V\right)\times \vc V_M\cdot\psi}\,.
\end{align*}

At this point, remark that, since $\vc V$ is a weak-limit point of the sequence $\big(\vc V_\eps\big)_\veps$, in view of Corollary \ref{c:momentum} and Lemma \ref{l:momentum}, it enjoys
the following property: 
\[
\vc V\,\in\,L^\infty_{\rm loc}\big(\R_+;L^2(\Omega)\big)\,\cap\,L^2_{\rm loc}\big(\R_+;W^{1,p_1}(\Omega)\big)\,,\qquad\qquad p_1\,:=\,\frac{2\k}{\k+1}\,.
\]
In particular, we also have $\vc V\,\in\,L^2_{\rm loc}\big(\R_+;L^{p_2}(\Omega)\big)$, where $p_2\,:=\,6\k/(\k+3)$.

Hence, repeating the computations used in the final part of the proof to Lemma \ref{l:curl}, we get that
\[
 \lim_{M\ra+\infty}\intTO{S_M\curl\left(\frac{1}{b}\,\vc V\right)\times \vc V_M\cdot\psi}=\,
 \lim_{M\ra+\infty}\intTO{\curl\left(\frac{1}{b}\,\vc V_M\right)\times \vc V_M\cdot\psi}\,,
\]
and performing computations in \eqref{eq:comp-conv} backwards, we deduce that, for any test function $\psi$ such that $\div(b\,\psi)=0$, we have
\[
 \lim_{M\ra+\infty}\intTO{S_M\curl\left(\frac{1}{b}\,\vc V\right)\times \vc V_M\cdot\psi}\,=\,
\lim_{M\ra+\infty}-\intTO{\frac{1}{b}\,\vc V_M\otimes \vc V_M:\nabla\psi}\,,
\]
where we have also used the fact that $\div \vc V\,=\,0$ (as it follows from taking the limit in the mass equation, recall Proposition \ref{p:constraints} above).

Now, using that $\vc V\in L^2_T(W^{1,p_1})$ and arguing as in \eqref{est:Id-S}, it is easy to see that, for almost any $t\in[0,T]$, one has
\[
 \left\|S_M\big(\vc V(t)\big)\,-\,\vc V(t)\right\|_{L^{p_1}}\,\lesssim\,2^{-M}\,\left\|\Grad\vc V(t)\right\|_{L^{p_1}}\,,
\]
which immediately implies that
\[
S_M\vc V\,\longrightarrow\,\vc V\qquad\qquad \mbox{ strongly in }\qquad L^2_T(L^{p_1})\,. 
\]
In fact, this convergence holds true even in $L^2_T(B^1_{p_1,2})$ (notice that, by \eqref{embedding} below, we have $W^{1,p_1}\hookrightarrow B^{1}_{p_1,2}$),
owing to Lemma \ref{l:Id-S} and  the Lebesgue dominated convergence Theorem; however, the previous weaker convergence result is enough for our scopes.

Therefore, $S_M\vc V$ converges strongly to $\vc V$ in any intermediate space between $L^2_T(L^{p_1})$ and $L^2_T(L^{p_2})$, thus also in $L^2_T(L^2)$ for instance.
Thanks to this latter property, we can compute
\[
 \lim_{M\ra+\infty}-\intTO{\frac{1}{b}\,\vc V_M\otimes \vc V_M:\nabla\psi}\,\dt\,=\,-\intTO{\frac{1}{b}\,\vc V\otimes \vc V:\nabla\psi}\,\dt\,.
\]

\subsection{Deriving the asymptotic system: final computations} \label{ss:final-comp}

In Subsections \ref{ss:approx} and \ref{ss:comp-comp}, we have seen how passing to the limit in the convective term and the pressure term.
About the latter, we recall that we have to make use of Lemma \ref{l:pressure}, and more precisely of the relation
\[
 \frac{1}{\veps^2}\,\Big(\nabla p(\vr_\eps)\,-\,\vr_\ep\nabla G\Big)\,=\,\frac{1}{\veps}\,b\,\nabla\big(H''(b)\,\phi_\ep\big)\,+\,\frac{1}{\veps^2}\nabla\Pi(\vr_\ep;b)\,,
\]
where the first term on the right disappears whenever tested again a test function satisfying \eqref{eq:test}, whereas the second term
is combined with the convective term to give rise to small remainders, in the sense of relations \eqref{eq:def-rem} and \eqref{eq:def-rem_int}.

On the other hand, the same computations performed
in Proposition \ref{p:constraints} show how dealing with the continuity equation for the densities $\vr_\veps$ and with the time derivative term $\big(\d_t(\vr_\veps\,\vu_\veps)\big)_\veps$ in
the momentum equation.
Therefore, in order to complete the proof of Theorem \ref{thm:main}, we must show convergence of the viscosity term
\[
\nu\intTO{\vre\,\tD\vue:\Grad\psi}\,,
\]
where $\psi$ is as in \eqref{eq:test} and is such that $\Supp\psi\subset[0,T[\,\times \Omega$.
We start by observing that only the integral over $\Omega^\ep_\ess$ matters, owing to the uniform bounds $\big(\sqrt{\vre}\,\Grad\vue\big)_\veps\,\subset\,L^2_T(L^2)$ and
$\big(\big[\sqrt{\vre}\big]_\res\big)_\veps\,\subset\,L^\infty_T(L^{2\g})$, with $1/2\,+\,1/(2\g)\,<\,1$.

Next, on $\Omega_\ess^\ep$ we use the strong convergence $\vre\to b$ in $L^\infty_T(L^2)\cap L^2_T(L^6)$ and the weak convergence $\tD\vue\,\rightharpoonup\,\tD \vc U$
in $L^2_T(L^{p_1})$. We observe that $1/6\,+\,1/p_1\leq 1$. Therefore, we deduce that, for any test function $\psi$ as above, we have
\[
\nu\intTO{\vre\,\tD\vue:\Grad\psi}\,\longrightarrow\,\intTO{b\,\tD\vc U:\Grad\psi}\qquad\qquad\mbox{ when }\qquad \veps\,\ra\,0^+\,.
\]

Theorem\ref{thm:main} is now  proven. \qed

%%%%%%%%%%%%%%%%%%%%%%%%%%%%%%%%%%%%%%%%%%%%%%%%%%%%%%%%%%%%%%%%%%%%%%%%%%%%%%%%%%%%%%%%%%%%%%%%%%%%%
\appendix

\section{Appendix: elements of Fourier analysis} \label{app:LP}

We recall here the main ideas of Littlewood-Paley theory, which we will exploit in our analysis.
The classical construction is usually given in the $\R^d$ setting: we refer \tsl{e.g.} to Chapter 2 of \cite{B-C-D} for details.
However, everything can be adapted (see \tsl{e.g.} reference \cite{Danchin}) to cover also the case of a $d$-dimensional periodic box $\T^d_a$, where $a\in\R^d$
(this means that the domain is periodic in space with, for any $1\leq j\leq d$, period equal to $2\pi a_j$ with respect to the $j$-th component).

For simplicity of presentation, we focus here on the case where $a_j=1$ for all $j$, and we simply write the spacial domain as $\T^d$. We also denote by
$|\T^d|\,=\,\mc L(\T^d)$ the Lebesgue measure of the box $\T^d$.

\medbreak
First of all, let us recall that, for a tempered distribution $u\in\mc S'(\T^d)$, we denote by $\mc Fu\,=\,\big(\what u_k\big)_{k\in\Z^d}$ its Fourier series,
so that we have
\[
u(x)\,=\,\frac{1}{|\T^d|^{1/2}}\,\sum_{k\in\Z^d}\what{u}_k\,e^{ik\cdot x}\,.
\]

Next, we introduce the so called \textit{Littlewood-Paley decomposition}, based on a non-homogeneous dyadic partition of unity with
respect to the Fourier variable. 
We fix a smooth scalar function $\vphi$ such that $0\leq \vphi\leq 1$, $\vphi$ is even and supported in the ring $\left\{r\in\R\,\big|\ 5/6\leq |r|\leq 12/5 \right\}$, and such that
\[
\forall\;r\in\R\setminus\{0\}\,,\qquad\qquad \sum_{j\in\Z}\vphi\big(2^{-j}\,r\big)\,=\,1\,.
\]

Let us define $|D|\,:=\,(-\Delta)^{1/2}$ as the Fourier multiplier\footnote{Throughout we agree  that  $f(D)$ stands for 
the pseudo-differential operator $u\mapsto\mc{F}^{-1}(f\,\mc{F}u)$.} of symbol $|k|$, for $k\in\Z^d$.
The dyadic blocks $(\Delta_j)_{j\in\Z}$ are then defined by
$$
\forall\;j\in\Z\,,\qquad\qquad
\Delta_ju\,:=\,\varphi(2^{-j}|D|)u\,=\,\sum_{k\in\Z^d}\vphi(2^{-j}|k|)\,\what u_k\,e^{ik\cdot x}\,.
$$
Notice that, for $j<0$ negative enough (in general, depending on the box $\T^d_a$), one has $\Delta_j\equiv0$. In addition, one has the following Littlewood-Paley decomposition in $\mc S'(\T^d)$:
%for any $u\in\mc{S}'(\T^d)$,  one has the equality 
$$
\forall\;u\in\mc{S}'(\T^d)\,,\qquad\qquad u\,=\,\what u_0\,+\, \sum_{j\in\Z}\Delta_ju\qquad\mbox{ in }\quad \mc{S}'(\T^d)\,.
$$

Finally, we introduce the following low frequency cut-off operators: for any $j\in\Z$, we define
\begin{equation} \label{eq:S_j}
S_ju\,:=\,\what u_0\,+\,\sum_{m\leq j-1}\Delta_{m}u\,.
\end{equation}

We explicitly remark that, for any $j\in\Z$, the operators $\Delta_j$ and $S_j$ are linear operators which are bounded on $L^p$ for any $p\in[1,+\infty]$, with norm \emph{independent} of $j$
and $p$.

At this point, we present a simplified version of the classical \emph{Bernstein inequalities}, which turns out to be enough for our scopes. We refer to Chapter 2 of \cite{B-C-D}
for the statement in its full generality.

\begin{lemma} \label{l:bern}
 There exists a constant $C>0$, only depending on the space dimension $d$, on the size of the torus $\T^d_a$ and on the support of the function $\vphi$ fixed above, such that the following properties
 hold true:
for any $j\in\Z$,
for any $\alpha\in\N^d$, for any couple $(p,q)\in[1,+\infty]^2$ such that $p\leq q$,  and for any smooth enough $u\in \mc S'(\T^d)$,  we  have
\begin{align*}
&\left\|\nabla^\alpha S_ju\right\|_{L^q}\, \leq\,
 C^{|\alpha|+1}\,2^{j|\alpha|+jd\left(\frac{1}{p}-\frac{1}{q}\right)}\,\|S_ju\|_{L^p} \\[1ex]
&\qquad\qquad\qquad\qquad\mbox{ and }\qquad\qquad\qquad
C^{-|\alpha|-1}\,2^{-j|\alpha|}\,\|\Delta_ju\|_{L^p}\,\leq\,
\|\nabla^\alpha \Delta_ju\|_{L^p}\,\leq\,C^{|\alpha|+1} \, 2^{j|\alpha|}\,\|\Delta_ju\|_{L^p}\,,
\end{align*}
where we have denoted $|\alpha|\,:=\,\sum_{j}\alpha_j$.
\end{lemma}   

By use of Littlewood-Paley decomposition, we can now define the class of Besov spaces.
\begin{df} \label{d:B}
Let $s\in\R$ and $1\leq p,r\leq+\infty$. The \emph{non-homogeneous Besov space} $B^{s}_{p,r}\,=\,B^s_{p,r}(\T^d)$ is the  set of tempered distributions $u\in\mc S'(\T^d)$ for which
$$
\|u\|_{B^{s}_{p,r}}\,:=\,\left(\left|\what u_0\right|^r\,+\,\sum_{j\in\Z}2^{jsr}\,\|\Delta_ju\|^r_{L^p}\right)^{1/r}\,<\,+\infty\,,
$$
with the standard modification in the definition of the norm in the case when $r=+\infty$.
\end{df}

It is well known that, for all $s\in\R$, the space $B^s_{2,2}$ coincides with $H^s$, with equivalent norms:
\begin{equation} \label{eq:LP-Sob}
\|f\|^2_{H^s}\,\sim\,\left|\what u_0\right|^2\,+\,\sum_{j\in\Z}2^{2sj}\,\|\Delta_ju\|_{L^2}^2\,.
\end{equation}
When $p\neq 2$, non-homogeneous Besov spaces are interpolation spaces between Sobolev spaces $W^{k,p}$:
for all $p\in\,]1,+\infty[\,$, one has the chain of following continuous embeddings:
\begin{equation} \label{embedding}
B^0_{p,\min(p,2)}\,\hookrightarrow\, L^p\,\hookrightarrow\, B^0_{p,\max(p,2)}\,.
\end{equation}

As an immediate consequence of the Bernstein inequalities, one gets the following Sobolev-type embedding result.
\begin{prop}\label{p:embed}
Let $1\leq p_1\leq p_2\leq+\infty.$ 
The, the space $B^{s_1}_{p_1,r_1}$ is continuously embedded in the space $B^{s_2}_{p_2,r_2}$ whenever
$$
s_2\,<\,s_1-d\left(\frac{1}{p_1}-\frac{1}{p_2}\right)\qquad\mbox{ or }\qquad
s_2\,=\,s_1-d\left(\frac{1}{p_1}-\frac{1}{p_2}\right)\;\;\mbox{ and }\;\;r_1\,\leq\,r_2\,. 
$$
\end{prop}

We conclude this appendix by recalling Lemma 2.73 of \cite{B-C-D}.
\begin{lemma} \label{l:Id-S}
If $1\leq r<+\infty$, for any $f\in B^s_{p,r}$ one has
$$
\lim_{j\ra+\infty}\left\|f\,-\,S_jf\right\|_{B^s_{p,r}}\,=\,0\,.
$$
\end{lemma}

\end{document}